\documentclass[12pt,a4paper]{article}
\usepackage{amsmath,pifont,amsfonts,amssymb,latexsym,multicol}
\usepackage{float,xspace,calc,theorem,ifthen}
\usepackage{fancyhdr,enumerate,psboxit,epsfig}
\usepackage{times,epic,amscd} 
\usepackage{version,xr}
\usepackage{aliascnt}
\usepackage{hyperref}
\usepackage{rotating}
\usepackage{color}

\hypersetup{colorlinks,pdfauthor={Andreas Knauf}}
\numberwithin{equation}{section}
\numberwithin{figure}{section}
\theoremstyle{change}
\newtheorem{theorem}{Theorem} [section]
\newtheorem {lemma}[theorem]{Lemma}
\newtheorem {prop}[theorem]{Proposition}

{\theorembodyfont{\normalfont}\newtheorem {remark}[theorem]{Remark}}
{\theorembodyfont{\normalfont}\newtheorem {remarks}[theorem]{Remarks}}
{\theorembodyfont{\normalfont}\newtheorem {example}[theorem]{Example}}
\newcommand{\beq}{\begin{equation}}
\newcommand{\eeq}{\end{equation}}
\newcommand{\Leq}[1]{\label{#1}\end{equation}}
\newcommand{\beqn}{\begin{eqnarray}}
\newcommand{\eeqn}{\end{eqnarray}}
\newcommand{\beqno}{\begin{eqnarray*}}
\newcommand{\eeqno}{\end{eqnarray*}}
\newcommand{\es}{\emptyset}
\renewcommand {\l}{\left}
\newcommand {\ri}{\right}
\newcommand {\vep}{\varepsilon}
\newcommand {\vv}{\varphi} 
  
\newcommand {\la}{\lambda}     

\newcommand {\LA}{\left\langle}
\newcommand {\RA}{\right\rangle}

\newcommand {\eh}{{\textstyle \frac{1}{2}}}
\newcommand {\ed}{{\textstyle \frac{1}{3}}}
\newcommand {\ev}{{\textstyle \frac{1}{4}}}
\renewcommand {\dh}{{\textstyle \frac{3}{2}}}
\newcommand {\ar}{\rightarrow}

\newcommand {\bC}{{\mathbb C}}

\newcommand {\bF}{{\mathbb F}}

\newcommand {\bN}{{\mathbb N}}
\newcommand {\bR}{{\mathbb R}}
\newcommand {\bZ}{{\mathbb Z}}

\newcommand {\bP}{{\mathbb P}}
\newcommand {\bQ}{{\mathbb Q}}
\newcommand{\rstr}{{\upharpoonright}}
\newcommand{\idty}{{\rm 1\mskip-4mu l}} 
\newcommand{\cB}{{\cal B}} 
\newcommand{\cC}{{\cal C}} %
\newcommand{\cD}{{\cal D}} %
\newcommand{\cF}{{\cal F}} %
\newcommand{\cH}{{\cal H}}

\newcommand{\cK}{{\cal K}}

\newcommand{\cN}{{\cal N}}

\newcommand{\cR}{{\cal R}} 

\newcommand{\bem}{\l(\! \begin{array}}
\newcommand{\eem}{\end{array}\!\ri)}
\newcommand{\bsm}{\left(\begin{smallmatrix}} 
\newcommand{\esm}{\end{smallmatrix}\right)}  

\newcommand{\qmbox}[1]{\quad\mbox{#1}\quad}
\renewcommand {\max}{{{\rm max}}}

\newcommand{\supp}{{\rm supp}}
\newcommand {\Zt}{\tilde{Z}}  

\newcommand{\hA}{{\bf A}}
\newcommand{\hB}{{\bf B}}
\newcommand{\Hi}{{\cal H}}
\newcommand{\hI}{{\bf I}}
\newcommand{\hJ}{{\bf J}}
\newcommand{\hL}{{\bf L}}
\newcommand{\hO}{{\bf O}}

\newcommand{\hR}{{\bf R}}
\newcommand{\hT}{{\bf T}}
\newcommand{\htt}{{\bf t}}
\newcommand{\hU}{{\bf U}}

\newcommand{\hY}{{\bf Y}}
\newcommand{\unity}{{\setlength{\unitlength}{1em}
                    \begin{picture}(0.5,0.8)
                    \put(0,0){$1$}
                    \put(0.38,0){\line(0,1){0.65}}
                    \end{picture}
                   }}
\begin{document}
\title {The Spectrum of an Adelic Markov Operator}
\author{Andreas Knauf\thanks{
Department of Mathematics,
Friedrich-Alexander-University Erlangen-Nuremberg,
Cauerstr.\ 11, D-91058 Erlangen, 
Germany, \texttt{knauf@mi.uni-erlangen.de}}}
\date{May 27, 2013}
\maketitle
\begin{abstract}
With the help of the representation of ${\rm SL}(2,\bZ)$ on the 
rank two free
module over the integer adeles, we define the transition operator  
of a Markov chain. The real component of its
spectrum exhibits a gap, whereas the non-real component
forms a circle of radius $1/\sqrt{2}$.
\end{abstract}
\tableofcontents
%
\section{Introduction and Setup}
In this article we consider unitary representations of ${\rm SL}(2,\bZ)$.
With $\hL$ and $\hR$ representing $\bsm1&1\\ 0&1\esm$ 
respectively $\bsm1&0\\ 1&1\esm$, we analyze the Markov operator 
\beq
\hT:=\eh(\hL+\hR).
\Leq{hT1}
As we indicate in this introductory section, its spectrum depends on the
representation chosen and controls various 
number theoretical equipartition rates. 

In Section \ref{sec:spa} we shall
then do the spectral analysis for relevant representations.

\subsubsection*{The Number-Theoretical Spin Chain}
%
The Euler product and Dirichlet series
\beq
Z(s):=\frac{\zeta(s-1)}{\zeta(s)}
=\prod_{p\in\bP}\frac{1-p^{-s}}{1-p^{1-s}}=\sum_{n=1}^\infty \vv(n)n^{-s}
\Leq{first:Z} 
(with Euler's $\vv$--function) converge in the half-plane $\Re(s)>2$.

On the discrete abelian group 
\[G:= \bigoplus_\bN (\bZ/2\bZ)\] 
with $\bZ/2\bZ\cong \big(\{0,1\},+\big)$ we define
$h: G\to \bN$ by $h(0):=1$ and
\beq
h(g_1,\ldots ,g_{n-1},1,0,\ldots) :=
h(g_1,\ldots ,g_{n-1},0,\ldots)+h(1-g_1,\ldots, 1-g_{n-1},0,\ldots).
\Leq{recurs}
Denoting by $h_k$ the restriction of $h$ to the subgroup 
$(\bZ/2\bZ)^k \cong 
\big(\!\bigoplus_{\ell=1}^k \bZ/2\bZ\big)\oplus 
\big(\!\bigoplus_{\ell=k+1}^\infty\{0\}\big)$ of $G$,
{\em e.g.} for $k=3$ one obtains
\begin{center}
\begin{tabular}{c|cccccccc}
   $g$ & 000 & 001 & 010 & 011 & 100 & 101 & 110 & 111\\ \hline
  $h_3(g)$ &  1 & 4 & 3 & 5 & 2 &5&3&4
 \end{tabular}\ .
\end{center}
The Dirichlet series equals
\[Z(s)=\sum_{g\in G} h(g)^{-s} = \sum_{g\in G} \exp(-s H(g)),\]
with $H:=\ln(h): G\to [0,\infty)$, see \cite{Kn1}. 
$Z$ is approximated by
\[Z_k(s):=
\sum_{g\in (\bZ/2\bZ)^k} h_k(g)^{-s}
= \sum_{g\in (\bZ/2\bZ)^k} \exp\big(\!-s H_k(g)\big) \qquad(s\in\bC;\ k\in\bN_0),\]
for $H_k:=\ln(h_k)$.
Writing $Z_k(s) = \sum_{n=1}^\infty \vv_k(n)n^{-s}$ defines the $\vv_k:\bN\to\bN_0$.

Since $\vv_k\le\vv$, but $Z$ has a pole at 2, $Z_k$ absolutely converges to $Z$
exactly in the half plane $\Re(s)>2$.

Although the Pontryagin dual $\widehat G$ of the direct sum $G$ is
the compact direct product $\prod_\bN (\bZ/2\bZ)$, it is impossible
to define the Fourier transform of $H$ in the sense of locally compact
abelian groups, since $H$ is not bounded, and $\ell^p(G)\subseteq \ell^\infty(G)$
for $1 \le p < \infty$.
However, the Fourier transform of $-H$ 
exists in the sense
\[j:G\setminus\{0\}\to\bR\qmbox{,}
j(t):= - \lim_{k\ar\infty} \  2^{-k}\!\!\!\! \sum_{g\in (\bZ/2\bZ)^k} 
H_k(g)\,(-1)^{\LA g,t\RA}\]
(in fact it can be defined in the same way on $\widehat G\setminus\{0\}$, 
but then vanishes outside $G\subseteq \widehat G$).
\begin{remark}[Statistical mechanics interpretation] 
In the language of statistical mechanics and
thermodynamic formalism (see, {\em e.g.} Ruelle \cite{Ru}), 
$Z$ is the {\em partition function} 
of an infinite spin chain with {\em energy function} $H$. Similarly, $Z_k$
is the partition function of $k$ spins.

The {\em interaction} $j$ has the following properties (see \cite{Kn1,Kn2,GK}):
\begin{itemize}
\item 
it is asymptotically translation invariant,
\item 
of long range (leading to phase transition),
\item 
and is non-negative.
\end{itemize}
The last property is called {\em ferromagnetism}.
Ferromagnetic spin systems are known to obey many specific inequalities.

For ferromagnetic {\em Ising} systems (that is, spin systems
that unlike our $j$ possesses only pair interaction), 
by the Lee-Young Theorem 
the partition function zeroes are on the imaginary line of the complex
plane of magnetization.~\hfill$\Diamond$
\end{remark}
Since $\vv(n)=|\{k\in\{1,\ldots,n\}\mid {\rm gcd}(k,n)=1\}|\ (n\in\bN)$,
we can consider the Dirichlet series $Z$ as a sum over 
$\Lambda:=\{\bsm a\\b\esm\in\bZ^2\mid{\rm gcd}(a,b)=1\}$:
\beq
Z(s)=1+\sum_{\mbox{\tiny $\bsm a\\ b\esm$}\;\in\;\Lambda\cap\bN^2} (a+b)^{-s}\qquad
\big(\Re(s)>2\big) .
\Leq{Z:2}
But with the unitary representation of ${\rm SL}(2,\bZ)$ on the Hilbert space 
$\ell^2(\Lambda)$ we may also write
\[Z(s)\!=\!1\!+\!\textstyle\sum_{k=0}^\infty\! 
\LA (2\hT_{\Lambda})^k\idty_{\mbox{\tiny $\bsm1\\1\esm$}},  \|\cdot\|_1^{-s}
\RA_{\!\ell^2(\Lambda)},\]
with $\hT$ generally defined in \eqref{hT1}. This first representation is
analyzed in Subsection \ref{subs:rc}.
\subsubsection*{Twisting the Dirichlet Series} 
%
%
For the Dirichlet series $Z$ the half-plane of convergence is not larger
than the half-plane $\{\Re(s)>2\}$ of absolute convergence. To have
a chance to look into the critical strip,
we now twist the Euler product of $Z$ in (\ref{first:Z}) to obtain a Dirichlet series 
\beq
\Zt(s) := \prod_{p\in\bP}\frac{1+p^{-s}}{1+p^{1-s}}
=\sum_{n=1}^\infty \lambda(n)\, \vv(n)\, n^{-s}
= \frac{\zeta(s)\,\zeta\big(2(s-1)\big)}{\zeta(s-1)\,\zeta(2s)}\ \quad\big(\Re(s)>2\big)
\Leq{Z:tilde}
with the Liouville function, given by 
$\lambda(p_1^{a_1}\cdot\ldots\cdot p_k^{a_k})=(-1)^{a_1+\ldots+a_k}$ for $p_i\in\bP$.
This has the following properties:
\begin{itemize}
\item 
Of the four zeta functions appearing in (\ref{Z:tilde}), only $s\mapsto \zeta(s-1)$
is not absolutely convergent for $\Re(s)>3/2$.
\item 
The pole of $\zeta$ at $s=1$, leading to the pole of
$Z$ at $s=2$, gives rise to $\Zt(2)=0$. 
\item 
The non-trivial zeros of $\zeta$, shifted by 1 for $Z$, now turn into poles of $\Zt$. 
\item
$\Zt$ has an additional pole at $3/2$.
\end{itemize}

So the Dirichlet series $\Zt$ converges in the half-plane 
$\{s\in\bC\mid\Re(s)>s_0+1\}$
if and only if there are no zeros of $\zeta$ with real part 
strictly larger than $s_0\ge 1/2$. 

Instead of considering convergence of the partial sums 
\beq
\hat{Z}_N(s) := \textstyle \sum_{n=1}^N \lambda(n)\, \vv(n)\, n^{-s}\qquad(s\in\bC),
\Leq{hat:Z}
to $\Zt(s)=\sum_{g\in G}\lambda\circ h(g) \ h(g)^{-s}$, one may also look at the convergence of 
\beq
\Zt_k(s) := \sum_{g\in (\bZ/2\bZ)^k} \lambda\circ h_k(g)\; h_k(g)^{-s}
=\textstyle \sum_{n=1}^\infty \lambda(n)\, \vv_k(n)\, n^{-s}.
\Leq{Zt1}
If we set $N\equiv N(k):=\lfloor \pi\sqrt{2^{k}/3}\rfloor$ in \eqref{hat:Z}, 
then both sums include asymptotically 
the same number of terms $n^{-s}$  (since $\sum_{n=1}^N \vv(n)\sim 3/\pi^2\ N^2$).

\begin{figure}[ht]
\begin{center}
\hspace{-2mm}

\includegraphics[height=139mm,angle=90]
{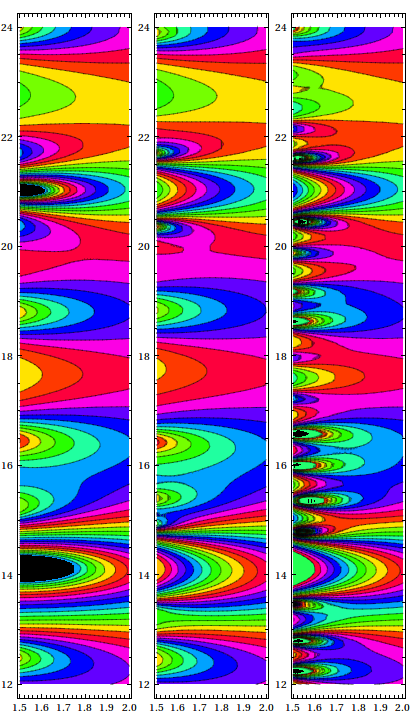}%
\caption{Bottom graph: Modulus of the 
twisted Dirichlet function $\Zt$ in its critical (half-) strip 
$\Re(s)\in[3/2,2]$. The real (respectively imaginary) part of $s$ is plotted
on the $y$--axis ($x$--axis).
Poles (shifted Riemann zeroes) of $\Zt$ are located
at the points $3/2+14.1\,\imath$ and $3/2+21.0\,\imath$ of the abscissa.\newline
The approximants $\Zt_k$ (middle)
respectively $\hat{Z}_{N(k)}$ (up) of $\Zt$, for $k=25$.}
\label{ZZZ}
\end{center}
\end{figure}

The numerics of their moduli (see Figure \ref{ZZZ}, and \cite{Kn3}) 
suggests that {\em if} there is convergence in the half-plane $\Re(s)>3/2$
(that being equivalent to the truth of RH), the convergence
of $\Zt_k$ may be even better than for $\hat{Z}_{N(k)}$. 

$\Zt_k$ numerically converges better, since
it has better self-averaging properties
than $\hat{Z}_{N}$ (which is of course a very non--smooth truncation 
of~$\Zt$).

A heuristic reason for such a supposed convergence of $\hat{Z}_N$ is to compare
the terms  $\lambda(n)$ appearing in $\hat{Z}_N$ to i.i.d.\ random 
variables which take the
values $\pm 1$ with equal probability $\eh$. 
For the case of $1/\zeta(s)=\sum_{n=1}^\infty \mu(n) n^{-s}$
a well-known similar heuristic goes back to Denjoy  (1931), and is described
in Section 12.3 of Edwards \cite{Ed}. 

Although that is obviously absurd in the literal sense,
we show in this article that there is some truth to the argument.
In the case $\Zt_k$ we thus consider the ensembles 
$g\mapsto\lambda\circ h(g)$ for
$g\in (\bZ/2\bZ)^k$. Let us begin with a simpler question.
\begin{example}[divisibility]
To convey the idea, we ask about the divisibility properties of the ensemble
$h_k(g) \quad (g\in (\bZ/2\bZ)^k)$ of $2^k$ integers. For division by $n=3$ a statistic
is given in Table~\ref{statist}. 
\begin{table}[ht]
\begin{center}
{\tiny 
\[
\begin{array}{cccccccccccccccccc}
\!\!\! 0 &\! 1 &\! 2 &\! 3 &\! 4 &\! 5 &\! 6 &\! 7 &\! 8 &\! 9 &\! 10 &\! 11 &\! 12 &\! 13 &\! 14 
&\! 15 &\! 16 &\! 17\!\!\! \!  \\
\!\!\!  0 &\! 0 &\! 2 &\! 2 &\! 2 &\! 10 &\! 18 &\! 26 &\! 66 &\! 138 &\! 242 &\! 506 
&\! 1058 &\! 2026 &\! 4050 &\! 8282 &\! 16386 &\! 32586\!\!\! \!  \\
\!\!\! -\frac{1}{4} &\! -\frac{1}{2} &\! 1 &\! 0 &\! -2 &\! 2 &\! 2 &\! -6 &\! 2 &\! 10 
&\! -14 &\! -6 &\! 34 &\! -22 &\! -46 &\! 90 &\! 2 &\! -182\!\!\! \!  \\
\!\!\! 1 &\! \sqrt{2} &\! 2 &\! 2 \sqrt{2} &\! 4 &\! 4 \sqrt{2} &\! 8 &\! 8 \sqrt{2} &\! 16 
&\! 16 \sqrt{2} &\! 32 &\! 32 \sqrt{2} &\! 64 &\! 64 \sqrt{2} &\! 128 &\! 128
   \sqrt{2} &\! 256 &\! 256 \sqrt{2}
\end{array}
\]
}
\caption{Statistic for divisibility by 3 of $h_k$. 
Upper row: exponent $k$ of the group cardinality $2^k$.
Second row: Number of $g\in (\bZ/2\bZ)^k$ with $3|h_k(g)$. 
Third row: After subtraction of expectation value $\ev 2^k$. 
Last row, for comparison: $2^{k/2}$, the square root of $|(\bZ/2\bZ)^k|$.}
\label{statist}
\end{center}
\end{table}
The set 
\[\Lambda(n):=\big\{\bsm a\\ b\esm\in(\bZ/n\bZ)^2\mid \gcd(a,b,n)=1\big\}
\qquad(n\in\bN)\] 
has the cardinality of Jordan's totient function
$J_2(n)=n^2\prod_{p\in \bP: p|n} (1-p^{-2})$, whereas its anti-diagonal 
$\{\bsm a\\ b\esm\in\Lambda(n)\mid a+b=0\}$ is of size $\vv(n)$.
So for the uniform distribution on $\Lambda(n)$ the expectation value of divisibility 
by $n$ equals the quotient $1/\psi(n)$ of these numbers, with Dedekind's psi
function $\psi$. For $n\in\bP$ this simplifies to $1/(n+1)$. 
Indeed, line 2 of Table~\ref{statist} is compatible
with the presumption that the quotient $1/4$ is approached as $k\ar\infty$.

Moreover, considering lines 3 and 4 of Table~\ref{statist}, 
the deviation from the expectation $|(\bZ/2\bZ)^k|/\psi(n)=2^{k-2}$ 
seems to be of the order $2^{k/2}$, the square root of the number of elements. 
This is similar to the sum
of $2^k$ i.i.d.\ random variables.

Both facts can be proven easily in the example at hand, using the matrix 
\beq \htt(3) :=
\eh\mbox{
{\tiny $\left(\begin{array}{ccccccccc}
 2 & 0 & 0 & 0 & 0 & 0 & 0 & 0 & 0 \\
 0 & 1 & 0 & 0 & 0 & 0 & 0 & 1 & 0 \\
 0 & 0 & 1 & 0 & 0 & 1 & 0 & 0 & 0 \\
 0 & 0 & 0 & 1 & 0 & 1 & 0 & 0 & 0 \\
 0 & 1 & 0 & 1 & 0 & 0 & 0 & 0 & 0 \\
 0 & 0 & 0 & 0 & 1 & 0 & 0 & 0 & 1 \\
 0 & 0 & 0 & 0 & 0 & 0 & 1 & 1 & 0 \\
 0 & 0 & 0 & 0 & 1 & 0 & 0 & 0 & 1 \\
 0 & 0 & 1 & 0 & 0 & 0 & 1 & 0 & 0
\end{array}\right)
$}}.
\Leq{mat3}
\begin{figure}
\centering
\includegraphics[width=0.49\columnwidth]{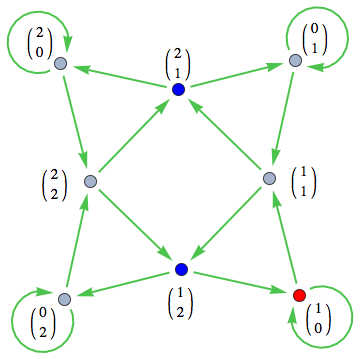}
\caption{The digraph with state space $\Lambda(3)$ 
for the Markov matrix $\htt(3)$.
The chain starts in the state ${1\choose 0}$. 
Divisibility by 3 of $h_k$ 
has a proportion equal to the sum of probabilities for the states 
${1\choose 2}$ and ${2\choose 1}$ after $k-1$ iterations of $\htt(3)$.}%
\label{fig:markov}%
\end{figure}
%
\begin{enumerate}[$\bullet$]
\item 
This acts on $\ell^2\big((\bZ/3\bZ)^2\big)$, is indexed lexicographically  
and is the Markov matrix for the recursion
(\ref{recurs}), taken (mod 3). 
So we obtain the second line of Table~\ref{statist} as matrix entries or, equivalently, 
$\ell^2$ scalar products:
\[ \big|\big\{g\in (\bZ/2\bZ)^k\mid h_k(g)=0\,(\mbox{mod }3)\big\}\big|
=\LA \big( 2\htt(3) \big)^k \,\idty_{\bsm 1\\0\esm},
\idty_{\bsm 1\\2\esm}+\idty_{\bsm 2\\1\esm}\RA.\]
\item 
The spectrum of $\htt(3)$ equals
\beq
\textstyle \left\{1,1,\frac{1}{4} \left(-1+i \sqrt{7}\right),
\frac{1}{4} \left(-1-i \sqrt{7}\right),\frac{1}{2},\frac{1}{2},\frac{1}{2},0,0\right\}.
\Leq{lambda3:spectrum}
The expectation $\frac{1}{4}2^k$ equals
$\LA\big( 2\htt(3)\big)^k P_{{\rm PF}}\, \idty_{\bsm 1\\ 0\esm},
\idty_{\bsm 2\\1\esm}+\idty_{\bsm 1\\2\esm}\RA$, with the orthogonal
projection $P_{{\rm PF}}$ to the Perron-Frobenius eigenspace 
${\rm span}(\idty_{\Lambda(3)})$.
\item 
On the other hand, apart from the two Perron-Frobenius
eigenvalues 1 in (\ref{lambda3:spectrum}), 
whose eigenfunctions are constant on the two parts of 
the $\htt(3)$--invariant decomposition 
$(\bZ/3\bZ)^2=\{\bsm 0\\0\esm\}\cup\Lambda(3)$ of the  state space, the spectral
radius of $\htt(3)$ is $1/\sqrt{2}$. As $\htt(3)$ is semisimple,
this leads to $\|(\htt(3)-P_{{\rm PF}})^k\|_{\Lambda(3)}\le c\,2^{-k/2}$, and
thus to the exponential estimate 
\[\Big|\big|\big\{g\in (\bZ/2\bZ)^k\mid h_k(g)=0\,(\mbox{mod }3)\big\}\big|-\ev 2^k\Big|
\ \le\ c \,2^{k/2}\qquad(k\in \bN)\]
for the deviation from the mean.
\hfill $\Diamond$
\end{enumerate}
\end{example}
Generalizing the example, the representations relevant for divisibility 
by $n\in\bN$ will be dealt with in Subsection \ref{subs:ZnZ}.

But the purpose of this article is to further generalize that kind of Markov
estimate, to gain some control on {\em joint} divisibility properties of the values
of $h_k$ ($k\in\bN$). Such a control is clearly necessary for
estimating the function $\lambda\circ h_k$ appearing in (\ref{Zt1}).
The natural language for this question is the one of adeles, and we consider
the corresponding representations in Subsections \ref{subs:fa} and \ref{subs:ac}.

Methods from the theory of expander graphs are used.
Some (like the one of Figure \ref{fig:markov}) but not all graphs arising here relate
to Ramanujan graphs, as defined by Lubotzky, Phillips,  and Sarnak in
\cite{LPS}. There is now a large literature on such expanders,
see, {\em e.g.}  \cite{BHV,DSV,HLW} for surveys.

\subsubsection*{Notation}
%
\begin{enumerate}[$\bullet$]
\item 
{\bf Miscellaneous:} The positive integers are $\bN=\{1,2,3\ldots\}$, whereas 
$\bN_0:=\{0,1,2\ldots\}$. With the primes $\bP:=\{2,3,5\ldots\}$ we set 
$\bP_\infty :=\{\infty\}\cup\bP$.

The imaginary unit is typeset as $\imath$. $\pm1$ is abbreviated by $\pm$.
Column vectors are sometimes written as row vectors, when this does no
harm.

We often apply a function to subsets $A$ of its domain: $f(A):=\{f(a)\mid a\in A\}$.

The spectrum of an operator $\bf O$ is denoted by ${\rm spec }({\bf O})$. So
its spectral radius $\sup (|{\rm spec }({\bf O})|)$ is bounded by 
its $L^2$ operator norm  $\|{\bf O}\|$.

We use the notation $\ell^2(M)$ for the $\bC$--Hilbert space 
over an infinite set $M$, 
equipped with counting measure, and similarly for the Banach spaces $\ell^p(M)$.
For finite sets $M$, in $\ell^p(M)$ the counting measure is often normalized to be a
probability measure.
\item 
{\bf Rings:} 
We index the absolute values 
$|\cdot |_v:\bQ_v\ar\bR$ by $v\in\bP_\infty$,
with $\bQ_\infty=\bR$. Whereas the ring of integers 
$\bZ_\infty=\bZ$ is  metrically discrete, for
$p\in \bP$ the ring $\bZ_p=\{x\in \bQ_p \mid |x|_p\le 1\}$ of 
$p$--adic integers is compact.

For a subset $S\subseteq \bP_\infty$ the ring $\bZ_S:=\prod_{v\in S}\bZ_v$
is equipped with the product topology.
$\bQ$ embeds diagonally in $\bZ_S$.  
For $S\subseteq \bP$ this embedding is
dense in the compact ring $\bZ_S$. This holds in particular for
$\bZ_\bP$ (also denoted by~$\widehat \bZ$).
\item 
{\bf Haar measures:}
The groups $\bZ_S$ are locally compact abelian and thus carry Haar measures:
$m_p$ on $(\bZ_p,+)$, normalized by $m_p(\bZ_p)=1\ (p\in\bP)$, 
$m_\infty$ (counting measure) on $(\bZ,+)$ and
\beq
m_S :=\bigotimes_{p\in S}m_p\qmbox{on} (\bZ_S,+). 
\Leq{def:mS}
So  $m_S(\bZ_S)=1$ if $S\subseteq \bP$.
\item 
{\bf Modules and Hilbert Spaces:} 
Our main objects will be the $\bZ_v$--modules $\bZ_v^2$, and 
$\bZ_S$--module $\bZ_S^2$ for $S\subseteq \bP_\infty$, 
with Haar measure $m_S^2:=m_S\otimes m_S$.
 
Since ${\rm SL}(2,\bZ)\subseteq {\rm SL}(2,\bZ_v)$ is a subgroup,
we obtain (via ${\bf O}_S f:= f\circ O^{-1}$) unitary representations 
of ${\rm SL}(2,\bZ)$ on the $\bC$--Hilbert spaces 
\[\cH_S:=L^2(\bZ_S^2,m_S^2)\qquad(S\subseteq \bP_\infty).\] 
\end{enumerate}
\subsubsection*{The Spectral Problem}
Given such a unitary representation on a Hilbert space $\cH$, rotation and inversion 
\[J,\, I:=J^2\;\in {\rm SL}(2,\bZ)\qmbox{,} 
J\bsm\ell\\ r\esm:=\bsm r\\-\ell\esm\ ,\quad I\bsm \ell\\ r\esm = \bsm-\ell\\-r\esm,\]
left and right addition
\[L, R\in {\rm SL}(2,\bZ)\qmbox{,} 
L\bsm \ell\\Êr\esm :=\bsm \ell+r \\ r\esm\ ,\quad 
R\bsm \ell\\ r\esm :=\bsm\ell \\ \ell+r\esm,\]
give rise to unitary operators 
$\hJ, \hI,\hL, \hR$ on $\cH$. 

Only when considering specific 
properties of a given unitary representation, we will provide 
the operators and spaces 
with a subindex, like $\hJ_S$ and $\cH_S$ for $S\subseteq \bP_\infty$.
We now temporarily omit that index for ease of notation.

We are to analyze the operators
\beq
\hT\in \hB(\cH)\qmbox{,} \hT= \eh (\hL + \hR).
\Leq{the:op}
Being convex combinations of unitaries, 
they are of operator norm $\| \hT \| \le 1$. 

Since we exclusively consider unitary representations of $O\in{\rm SL}(2,\bZ)$
given by measure preserving actions of the group via ${\bf O} f= f\circ O^{-1}$,
$\hT$ shares with $\hL$ and $\hR$ (and all above ${\bf O}$)
the invariance of ${\rm spec }(\hT)$ under complex conjugation.
%
\section{Spectral Analysis of the Operators $\hT$}\label{sec:spa}
%
We now analyze the operators (\ref{the:op}) on Hilbert spaces $\cH$ of 
unitary representations, more precisely on subspaces $\cH^\pm$ of $\cH$, 
defined by parity under inversion. 

In Section \ref{subs:ar} we find regions of
$\bC$ which are independent of the representation of ${\rm SL}(2,\bZ)$
and contain the spectrum (Proposition \ref{prop:algebraic}).
On the (more important) subspace $\cH^+$  we relate
the non-normal operator $\hT$ to a conjugate pair of projections,
thus preparing for a more precise spectral analysis
(Prop.\ \ref{prop:ABT}).

The rest of the section is about concrete representations:\\
For the regular representation (Sect.\ \ref{subs:rr}) the operator 
has a (minimal possible) spectral radius $1/\sqrt{2}$ (Prop.\ \ref{prop:SL}). 
This is shown using 
the Laplacians on a regular tree, and it leads to a similar statement
for the defining representation of  ${\rm SL}(2,\bZ)$
(Prop.\ \ref{prop:real} of Sect.\ \ref{subs:rc}).

To prepare for the finite adeles, Section \ref{subs:ZnZ} considers
representation of the finite groups ${\rm SL}(2,\bZ/n\bZ)$, $n\in\bN$.
Here a result by Bourgain and Varj\'u \cite{BV} of 2012 
is used to show the existence of an $n$-independent spectral gap.
Although $\hT$ is non-normal, the spectrum of the operator on the
finite adeles is the closure of the union of these spectra (see Prop.\ \ref{nonreal:ok}
of Sect.\  \ref{subs:ar}).
The final Section \ref{subs:ac} then combines these results and states that
the adelic operator $\hT$ has a gap.

%
\subsection{Algebraic Restrictions}\label{subs:ar}
%
In this subsection we continue to temporarily omit the subindex of the representation.
Since $I=\bsm -1&0\\ 0&-1 \esm\in {\rm SL}(2,\bZ)$ is a nontrivial involution, 
${\rm spec }(\hI)\subseteq\{-1,1\}$ for the spectrum 
of a unitary representation $\hI$. We thus obtain an orthogonal decomposition 
\beq
\cH=\cH^+\oplus \cH^-
\Leq{splitting}
of the Hilbert space into eigenspaces of $\hI$. 
Since $I$ is in the center of ${\rm SL}(2,\bZ)$,
the representations ${\bf O}$ of all $O\in{\rm SL}(2,\bZ)$
restrict to operators ${\bf O}^\pm$ on $\cH^\pm$. \\
In this article we mainly analyze the operators $\hT^+$. 

We first derive algebraic identities between operators ${\bf O}^\pm$ 
on $\cH^\pm$, that impose restrictions on ${\rm spec} (\hT^\pm)$ and
are valid for all unitary representations of ${\rm SL}(2,\bZ)$.

\begin{figure} 
\begin{center}
\begin{picture}(200,180)(-100,-90)
\put(-5,80){$\bsm\ell\\ r\esm$}
\put(45,40){$\bsm\ell+r\\ r\esm$}
\put(45,-40){$\bsm\ell+r\\ -\ell\esm$}
\put(-70,-40){$\bsm r\\ r-\ell\esm$}
\put(-70,40){$\bsm \ell\\ r-\ell\esm$}
\put(-5,-80){$\bsm r\\-\ell \esm$}
\put(15,75){\vector(2,-1){35}}
\put(-5,-65){\vector(-2,1){35}}
\put(-40,57){\vector(2,1){35}}
\put(50,-47){\vector(-2,-1){35}}
\put(65,25){\vector(0,-1){50}}
\put(-55,-25){\vector(0,1){50}}
\put(5,65){\vector(0,-1){130}}
\put(40,0){$R^{-1}$}
\put(-50,0){$L^{-1}$}
\put(40,65){$L$}
\put(-42,-65){$R$}
\put(-42,65){$R$}
\put(40,-65){$L$}
\put(10,0){$J$}
\end{picture}
\end{center}
\caption{Relations between the unitaries of left/right addition, and inversion}
\label{fig:hexagon}
\end{figure}
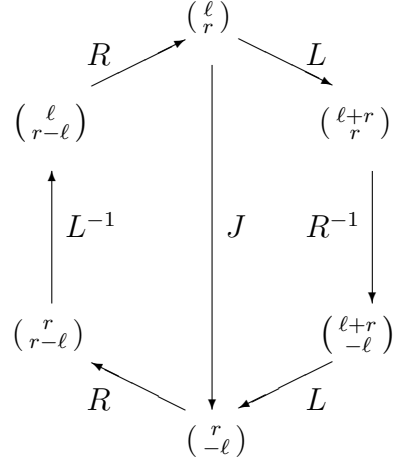

The restrictions of the representation $\hJ$ of
$J=\bsm 0&1\\ -1&0 \esm\in {\rm SL}(2,\bZ)$
to $\cH^+$ ($\cH^-$) are (anti-)selfadjoint, with
${\rm spec}(\hJ^+)\subseteq\{-1,1\}$ (and 
${\rm spec}(\hJ^-)\subseteq\{-\imath,\imath\}$).
We use the identities 
\beq
\hJ = \hL\hR^{-1}\hL = \hR^{-1}\hL\hR^{-1}
\Leq{le:hexagon}
for  the representations of
$L=\bsm 1&1\\ 0&1 \esm$ and $R=\bsm 1&0\\ 1&1 \esm$, 
exhibited in Figure \ref{fig:hexagon}.
So 
\beq
\hJ^* \hL\hJ= \hR^{-1}\qmbox{,}\hJ^* \hR\hJ= \hL^{-1}\qmbox{and thus}
\hT^*=\hJ^*\hT\hJ  \qmbox{on} \cH. 
\Leq{identi}
The equations (\ref{le:hexagon}) are complemented by
\beq
\hJ =\hI  \hL^{-1}\hR\hL^{-1} =\hI \hR\hL^{-1}\hR.
\Leq{le:hexagon2}
In order to proceed in the spectral analysis of $\hT$,
we introduce the operators 
\beq
\hY_+ := \ed \big(\dh \unity -\eh \hI+ \hR^{-1} \hL + \hL^{-1}\hR\big) 
\qmbox{and}
\hY_- := \hJ^*\hY_+\hJ
\Leq{Yplusmin}
on $\cH$. Note that the signs of $\hY_\pm$ do not refer to the splitting 
(\ref{splitting}), but instead 
\[\hY_{\pm}=\hY_{\pm}^+\oplus\hY_{\pm}^-\qquad 
\mbox{on }\cH=\cH^+\oplus \cH^-.\]
\begin{lemma}  \label{lem:algebraic}
$\hY_\pm$ are orthogonal projections, and
$\hT^*\hT = \frac{1}{8}(\idty+\hI+6\hY_+)$.
\end{lemma}
{\bf Proof.}
$\hY_\pm$ are self-adjoint, since $ \hL$ and $\hR$ are unitary and 
$\hI$ is an involution.\\
To show that $\hY_\pm^2=\hY_\pm$, one uses the relation 
$\hL\hR^{-1} \hL= \hI \hR\hL^{-1} \hR$, see Fig.~\ref{fig:hexagon}.
The second claim follows from
\[8\hT^*\hT=2(\hL^{-1}+\hR^{-1})(\hL+\hR) = 4\unity + 2(\hR^{-1}\hL+\hL^{-1}\hR)
= \idty+\hI+6\hY_+.\tag*{$\Box$}\]
%
%

\noindent
For the self-adjoint operators $\hA:=\hY_+-\hY_-$ and $\hB:= \idty-\hY_+-\hY_-$
(see Avron, Seiler and Simon \cite{ASS})
\beq
\hA^2+\hB^2=\idty\qmbox{,} \hA\hB + \hB\hA=0\qmbox{and}
[\hA^2,\hY_\pm]=[\hB^2,\hY_\pm]=0 .
\Leq{AB}
Like $\hY_\pm$, $\hA$ and $\hB$ are complexifications of real self--adjoint 
operators, so that, {\em e.g.} eigenfunctions may be chosen to be real.

Lemma \ref{lem:algebraic} implies that the kernels of 
$\hT^-$ and of $\hY_+^-$ coincide.
A more important use of the lemma is to localize the spectrum of $\hT^\pm$:
\begin{prop} \label{prop:algebraic} %
For all unitary representations of $\,{\rm SL}(2,\bZ)$ 
\begin{enumerate}
\item
$\hT^+$ is invertible, with $\ (\hT^+)^{-1} = 2(\hT^+)^* - \hJ^+$, 
and $2\hT^+ + (\hT^+)^{-1}$ is self-adjoint:
\beq
2\hT^+ + (\hT^+)^{-1} = -3\hJ^+\hB^+ = -3\hB^+\hJ^+.
\Leq{TJB:formula} 
So with the circle $\,C:= \big\{c\in \bC\mid |c|=1/\sqrt{2} \big\}$ 
and $I:=[-1,-\eh]\cup [\eh,1]$\\
\[{\rm spec}\big(\hT^+\big)\ \subseteq C\cup I    .\]
\item
$(\hT^-)^*\hT^-= \frac{3}{4}\hY_{+}^-$ so that with the disk 
$D := \big\{c\in \bC\mid |c|\le \sqrt{3}/2\big\}$
\[{\rm spec}\big(\hT^-\big) \subseteq D.\]
\end{enumerate}
\end{prop}
{\bf Proof.} 
For the sake of brevity we partly omit the superscript $\pm$.
\begin{enumerate}[1.]
\item
$\bullet$
We show that $\hT (2\hT^* - \hJ) = \unity$ on $\cH^+$, 
the identity $(2\hT^* - \hJ)\hT = \unity$
being similar: 
$2\hT\hT^*=\unity + \eh(\hL\hR^{-1}+\hR\hL^{-1})$,
and, using (\ref{le:hexagon}) and (\ref{le:hexagon2}) on $\cH^+$
\[2\hT \hJ = \hL\hJ +\hR\hJ = \hL(\hL^{-1}\hR\hL^{-1}) + \hR(\hR^{-1}\hL\hR^{-1})
= \hR\hL^{-1} + \hL\hR^{-1}.\]
$\bullet$
This implies (\ref{TJB:formula}), since its left hand side equals 
 \[ 2\hT^+ + (\hT^+)^{-1} = 2(\hT^++(\hT^+)^*) - \hJ^+ 
 = \hL^++ \hR^+ + (\hL^+)^{-1}+ (\hR^+)^{-1} - \hJ^+\]
whereas, reusing (\ref{le:hexagon}) and (\ref{le:hexagon2}),
$-3\hJ\hB=\hJ(-\hI + \hR^{-1}\hL + \hL^{-1}\hR + \hL\hR^{-1} + \hR\hL^{-1})
= \hI\hR+\hL+\hI\hL^{-1}+\hR^{-1}- \hI\hJ=-3\hB\hJ $ on $\cH^+$ 
since $\hI^+=\idty_{\cH^+}$.\\
$\bullet$
Since the preimage of $[-3,3]$ for the map $t\mapsto 2t+t^{-1}$ equals $C\cup I$,
we thus have to show that the resolvent set of $2\hT^++(\hT^+)^{-1}$ contains
$\bC\setminus[ -3,3]$. This follows  from (\ref{TJB:formula})
since $\|\hB\hJ\|\le \|\hB\| \|\hJ\| \le 1$, as $\|\hB\|^2=\|\hB^2\|\le 1$  by 
the first identity in (\ref{AB}).
\item
With Lemma \ref{lem:algebraic} we get $(\hT^-)^*\hT^-= \frac{3}{4}\hY_{+}^-$.
As $\| \hT^-\|^2 = \|(\hT^-)^*\hT^-\| =  \frac{3}{4}\|\hY_{+}^-\|\le  \frac{3}{4}$, 
the spectral radius of $\hT^-$ is smaller than $\sqrt{3}/2$.
\hfill$\Box$
\end{enumerate}
It might be instructive to have a look at the concrete spectra of $\hT^\pm$ for
a certain representation, depicted in Figure \ref{fig:spec} on Page \pageref{fig:spec}.

The two projections $\hY_\pm$ are unitarily equivalent and 
commute on the intersections of their kernels and ranges. 
We thus obtain (see Halmos in \cite{Ha2})
an orthogonal decomposition of the Hilbert space as a direct sum 
\[\cH= \cC_{0,0}\oplus\cC_{0,1}\oplus\cC_{1,0}\oplus\cC_{1,1}\oplus\widetilde\cN ,\]
with the closed {\em commuting} 
subspaces $\cC_{j,k}:= {\rm ker}(\hY_+- j\idty)\cap{\rm ker}(\hY_-- k\idty)$.

As we are analyzing $\hT$ instead of $\hY_\pm$, we modify that decomposition.
 From the formulae $\hT^*\hT = \frac{1}{8}(\idty+\hI+6\hY_+)$ 
(Lemma \ref{lem:algebraic}) and 
$\hT\hT^* = \frac{1}{8}(\idty+\hI+6\hY_-)$ we see that $\hT$  is 
normal on the domains $\cC_{0,0}$ and $\cC_{1,1}$, but not on
$\cC_{0,1}$ and $\cC_{1,0}$. Thus we set 
\beq
\cK:=  \cC_{0,0}\mbox{ , }\cR:=  \cC_{1,1}\qmbox{and obtain} 
\cH=\cK\oplus\cR\oplus\cN
\Leq{KRN}
in the sense of an orthogonal direct sum,
with $\cN= \cC_{0,1}\oplus\cC_{1,0} \oplus\widetilde\cN $.

In Example \ref{ex:not:semisimple} we exhibit a representation for which
$\hT^-\rstr_{\cN^-}$ is not semisimple (with $\cN^\pm:=\cN\cap \cH^\pm$).
As noted in Remark \ref{rem:jordan}.2 below, 
$\hT^+\rstr_{\cN^+}$ need not be semisimple, either.

Indeed, the spectral analysis of $\hT^+$ on the 
joint kernel respectively range subspaces $\cK^+=\cK\cap \cH^+$ and 
$\cR^+=\cR\cap \cH^+$ of $\hY_\pm^+$ is simple:
\begin{lemma}  \label{lem:algebraic2}
$\bullet$
$\hT^+=\eh (3 \hY_-^+ - \idty_{\cH^+}) \hJ^+ 
=\eh \hJ^+(3 \hY_+^+ - \idty_{\cH^+}) $.\\
$\bullet$
The operators
$\hJ^+$ and $\hT^+$ leave the decomposition 
$ \cH^+=\cK^+\oplus\cR^+\oplus\cN^+$ 
invariant.\\
$\bullet$
The operator $\hT^+$ restricts to the mutually orthogonal eigenspaces 
\[{\cal K}^{+,\pm}:=\{v\in {\cal K}^+ \mid \hJ^+\,v=\pm v\}\qmbox{and} 
{\cal R}^{+,\pm}:=\{v\in {\cal R}^+ \mid \hJ^+\,v=\pm v\}\] 
with $\,\hT^+\rstr_{\cK^{+,\pm}} = \mp\eh\idty\rstr_{\cK^{+,\pm}}$ and 
$\,\hT^+\rstr_{\cR^{+,\pm}} = \pm\idty\rstr_{\cR^{+,\pm}}$. \\
$\bullet$
Furthermore $\hT\rstr_{\cK^-}=\hT^*\rstr_{\cK^-}=0$.
\end{lemma}
{\bf Proof.} 
$\bullet$
On $\cH^+$ we have
$3 \hY_+ - \idty = \hR^{-1}\hL+ \hL^{-1}\hR$. 
Thus by (\ref{le:hexagon}) and (\ref{le:hexagon2}) 
$\hJ^+(3 \hY_+^+ - \idty)= \hI^+\hR^++\hL^+= 2\hT^+$, using  
$\hI^+=\idty_{\cH^+}$.\\
The definition (\ref{Yplusmin}) of $ \hY_-$ then implies
the other identity $\hT^+=\eh (3 \hY_-^+ - \idty_{\cH^+}) \hJ^+\!\!.$\\
$\bullet$
We can present $\cK$ as ${\rm ker}(\hY_+ + \hY_-)$ and $\cR$ as 
${\rm ker}(2\idty - \hY_+ - \hY_-)$. As $\hY_+ + \hY_-$ 
commutes with $\hJ$, we see that
both subspaces are  $\hJ$--invariant. Since $\hJ$ is 
self-adjoint, this also follows for the
orthogonal complement $\cN$ of $\cK\oplus\cR$.
$\hT^+$--invariance of $\cK^+$ and $\cR^+$ then follows from the formula
$\hT^+=\eh (3 \hY_-^+ - \idty_{\cH^+}) \hJ^+$.\\[1mm]
$\bullet$
With the same formula
$\hT^+\rstr_{\cK^{+,\pm}}= \eh (3 \hY_-^+ - \idty) \hJ^+\rstr_{\cK^{+,\pm}}
= \mp\eh\idty\rstr_{\cK^{+,\pm}}$ and
$\hT^+\rstr_{\cR^{+,\pm}}= \eh (3 \hY_-^+ - \idty) \hJ^+\rstr_{\cR^{+,\pm}}
= \pm\idty\rstr_{\cR^{+,\pm}}$.\\
$\bullet$
Finally, since $\hI^-=-\idty_{\cH^-}$, the formulae 
$\hT^*\hT = \frac{1}{8}(\idty+\hI+6\hY_+)$ and 
$\hT\hT^* = \frac{1}{8}(\idty+\hI+6\hY_-)$ imply 
$\hT^*\hT\rstr_{\cK^-} = \hT\hT^*\rstr_{\cK^-} = 0$ 
so that 
$\hT\rstr_{\cK^-}=\hT^*\rstr_{\cK^-}=0$.
\hfill$\Box$\\[2mm]

To analyze the rest of the spectrum one works on the last direct summand in
(\ref{KRN}),
the {\em non-commuting} subspace $\cN$.
Note again that on $\cN^+$ 
the operators $\hT^+$ and $\hJ^+$ do {\em not} commute, and hence 
$\hT^+\rstr_{\cN^+}$ is not normal. However, the deviation from normality 
can be controlled. 

To that aim we  note that, as $[\hB,\hJ]=0$, the 
self-adjoint operator $\hB^+$ decomposes as 
\beq
\hB^+=\hB^{+,+}\oplus \hB^{+,-}\qmbox{on}\cH^+= \cH^{+,+}\oplus\cH^{+,-} 
\Leq{bppm} 
for the subspaces
\[\cH^{+,\pm}:= \{v\in \cH^+ \mid \hJ^+\,v=\pm v\}.\]
Since $\hA^+\hJ^+ + \hJ^+\hA^+=0$, $\hA^+$ 
maps $\cH^{+,\pm}$ into $\cH^{+,\mp}$.
Actually, the restrictions 
\[\hA^+\rstr_{\cN^{+,\pm}}:\cN^{+,\pm}\ar \cN^{+,\mp}\]
for the orthogonal direct sum 
\beq
\hspace*{-2mm}\cN^+= \cN^{+,+}\oplus\cN^{+,-} 
\ \mbox{ with the subspaces }
 \cN^{+,\pm}:= \{v\in \cN^+ \mid \hJ^+\,v=\pm v\}
\Leq{splitting1}
are injective, since $\hA v=0$ implies $\hB^2 v=v$, using (\ref{AB}).

In the polar decomposition $\hO=\hU\, |\hO|$ of a {\em self-adjoint}
operator $\hO$, the partial isometry $\hU$ commutes with $|\hO|$, 
and is also called the {\em sign of} $\hO$ 
and written as ${\rm sign}(\hO)$, since ${\rm spec}(U)\subseteq \{-1,0,1\}$.

We now use the operator-valued spectral measure 
\[E:=E(\hO): \cB\big({\rm spec}(\hO)\big)\longrightarrow \hB(\cH)\]
of a normal operator $\hO\in\hB(\cH)$, defined on the Borel sigma algebra 
of its spectrum, whose unique existence is assured by the spectral theorem,
see {\em e.g.} Chapter IX, \S2 of Conway \cite{Co}. 
Since $\hB^{+}$ and $\hJ^{+}$ commute, we can use instead of  
$E(\hB^{+})$ its refinement
$E(\hB^{+,+})\oplus E(\hB^{+,-})$ (see Exercise 17 on p.\ 273 of \cite{Co}).
\begin{prop}\label{prop:ABT}
$\bullet$
$\hA^+$ transforms the spectral measures $E=E(\hB^{+,\pm})$ of the operators
$\hB^{+,\pm}$ as $\ \hA^+\;E(B)=E(-B)\;\hA^+\quad 
\big(B\in\cB({\rm spec}(\hB^{+,\pm}))\big)$.\\
$\bullet$
The operator $\hA^{+}\rstr_{\cN^+}$ splits orthogonally into
\beq
\hA^{+}=\hA^{+}_+\oplus \hA^{+}_-\qmbox{on} \cN^+=\cN^{+}_{+}\oplus\cN^{+}_{-}
\Leq{splitting2}
with $\cN^{+}_{\pm}:=P_\pm
(\cN^+)$ for the family ($P_+ + P_-=\idty_{ \cN^+}$) of orthogonal projectors 
\beq 
P_\pm :=\eh\big(\idty \pm {\rm sign}(\hA^{+}\rstr_{\cN^+})\big)\qmbox{on} \cN^+.
\Leq{AN:splitting}
$\bullet$
The spectral measure of $\hA^{+}_\pm$ on a measurable set
$\,0\not \in A\in \cB\big({\rm spec}(\hA^{+}_\pm)\big)$ equals
\beq
E(\hA^{+}_\pm)(A)= P_\pm\,E(|\hB^{+}|)
\textstyle{\big(\sqrt{1-A_\pm^2}}\,\big)\mbox{, with } 
A_\pm:=\{x\in A\mid \pm \,x\ge0\} .
\Leq{BA}
$\bullet$
Denoting by $b\mapsto dE^{+,\pm}(b)$ integration w.r.t.\ the spectral measure of 
$\hB^{+,\pm}$,
\beq
\hT^+\rstr_{\cN^+}= \int_{(-1,1)} \ev
\bsm 1-3b&3\sqrt{1-b^2}\\ -3\sqrt{1-b^2}& -1-3b\esm dE^{+,+}(b)\oplus dE^{+,-}(-b).
\Leq{rep:T}
This relates $\{(b,\pm1)\}\subseteq {\rm spec}(\hB^{+}, \hJ^{+})$ to 
$\big\{\ev\big(-3b\pm\sqrt{9b^2-8}\big)\big\} \subseteq  {\rm spec}(\hT^{+})$.

In particular, possibly excluding the values $\pm 1$ and $\pm \eh$, 
${\rm spec}(\hT^{+})$ is symmetric w.r.t.\ inversion at the circle 
$C$ of radius $1/\sqrt{2}$, including multiplicities.
\end{prop}
{\bf Proof.}
$\bullet$
The formula $\hA^+\,E(B)=E(-B)\,\hA^+$ follows with (\ref{AB}) from 
\[\hA^+f(\hB^+)=\hA^+\big(f_+(\hB^+)+f_-(\hB^+)\big) =
\big(f_+(\hB^+)-f_-(\hB^+)\big)\hA^+\]
for $f:\bR\ar\bR$ Borel measurable and $f_\pm(x):=\eh(f(x)\pm f(-x))$, 
by choosing $f:=\idty_B$
and noting that $(\idty_B)_+-(\idty_B)_- =\idty_{-B}$.\\
$\bullet$
(\ref {AN:splitting}) is an orthogonal projector family, 
since $P_++P_-=\idty_{\cN^+}$ and 
$\hA^{+}\rstr_{\cN^+}$  is injective, so that
\[\big({\rm sign}(\hA^{+}\rstr_{\cN^+})\big)^2=\idty_{\cN^+}\qmbox{and thus}
P_\pm^2 = P_\pm \mbox{ and }  P_\pm P_\mp = 0.\]
$\bullet$
Eq.\ (\ref{BA}) can be verified by merely using the relations 
$(\hA^{+})^2+(\hB^{+})^2=\idty$ and $\hA^{+}\hB^{+}+\hB^{+}\hA^{+}=0$ 
on $\cN^+$:\\
$|\hB^{+}| = \sqrt{(\hB^{+})^2}$ implies $[ \hA^{+},\, |\hB^{+}|\,]
=0=[ P_\pm,\, |\hB^{+}|\,]$. 
So $P_\pm\,E(|\hB^{+}|)(B)$ is self-adjoint for 
any Borel set $B$ and respects the decomposition given by (\ref {AN:splitting}).\\
Since  $P_\pm$ and $E(|\hB^{+}|)(B)$ are orthogonal projections, 
the right hand side of (\ref{BA}), applied to $B$, 
is a projection.
By the above, since $E(|\hB^{+,\pm}|)$ is a spectral measure, 
$P_\pm\,E(|\hB^{+}|)$ are spectral measures, too.

To verify (\ref{BA}), we substitute $|\hB^{+}|=\sqrt{\idty-(\hA^{+})^2}$.
Since $x\mapsto \sqrt{1-x^2}$ is strictly monotone on $[0,1]$ and on $[-1,0]$, 
$E(|\hB^{+}|)\textstyle{\big(\sqrt{1-A_\pm^2}}\,\big)=E(\hA^{+}_\pm)(A)$
on $\cN^+_\pm$.\\
$\bullet$
The representation (\ref{rep:T}) of $\hT^+$ based on the 
spectral measure of $\hB^{+,\pm}$ follows from
\beq
\hT^+=\ev(\idty-3\hA^+-3\hB^+)\hJ^+ = \ev\hJ^+(\idty+3\hA^+-3\hB^+)
\Leq{TAB}
(which is a consequence of the formula
$\hT^+=\eh (3 \hY_-^+ - \idty) \hJ^+ =\eh \hJ^+(3 \hY_+^+ - \idty) $
in Lemma \ref{lem:algebraic2}). The term $\ev(\idty-3\hB^+)\hJ^+$ of (\ref{TAB}) 
is diagonal 
in the common spectral resolution of $(\hB^+,\hJ^+)$ and thus 
gives rise to the diagonal of 
$ \ev\bsm 1-3b&3\sqrt{1-b^2}\\ -3\sqrt{1-b^2}& -1-3b\esm$.
The operator $\hA^+$ interchanges the two subspaces, and 
$(\hA^+)^2=\idty-(\hB^+)^2$ by (\ref{AB}), hence the off-diagonal entries.
\hfill $\Box$
\begin{remarks}\label{rem:jordan}
\begin{enumerate}[1.]
\item 
The splitting
(\ref{splitting1}) is diagonal  w.r.t.\ 
$\cN^+= \cN^{+}_+\oplus\cN^{+}_-$, 
whereas the splitting (\ref{splitting2}) is diagonal w.r.t.\ 
$\cN^+= \cN^{+,+}\oplus\cN^{+,-}$.
\item 
For the spectral parameter $b\in(-1,1)\setminus\big\{\pm\frac{\sqrt{8}}{3}\big\}$
the matrix $ \ev\bsm 1-3b&3\sqrt{1-b^2}\\ -3\sqrt{1-b^2}& -1-3b\esm$
appearing in (\ref{rep:T}) has eigenvalues 
$\frac{1}{4} \left(-3 b\pm \sqrt{9 b^2-8}\right)$ with eigenvectors
$\bsm\frac{-1\mp \sqrt{9 b^2-8}}{3 \sqrt{1-b^2}}\\1 \esm$.
We note that the vector entries $\frac{-1\mp \sqrt{9 b^2-8}}{3 \sqrt{1-b^2}}$ 
are of modulus one in the interval $|b| < \sqrt{8}/3$. 
For $b=\pm\frac{\sqrt{8}}{3}$ the degenerate eigenvalues $\pm1/\sqrt{2}$
are the points of the circle $C$ on the real axis, and
Jordan blocks arise.
\hfill $\Diamond$
\end{enumerate}
\end{remarks}
%
\subsection{The Regular Representation and Laplacians on Trees} \label{subs:rr}
%
Before we come to the analytic estimates of $\hT_S$ for sets 
$S\subseteq\bP_\infty$ of places,
we consider the corresponding problem for the operator $\hT_{\rm SL}$
on the Hilbert space 
\[\cH_{\rm SL}:= \ell^2\big({\rm SL}(2,\bZ)\big)\] 
with counting measure.
$\hT_{\rm SL}$ is defined by the unitary left regular representation
of ${\rm SL}(2,\bZ)$.
Following the notation in (\ref{splitting}) and (\ref{KRN})
we split into the orthogonal  subspaces 
\beq
\cH_{\rm SL} = \cH_{\rm SL}^+ \oplus \cH_{\rm SL}^-
\qmbox{and}
\cH_{\rm SL}^+=\cK_{\rm SL}^+\oplus\cR_{\rm SL}^+\oplus\cN_{\rm SL}^+.
\Leq{KRN:SL}

$\hT_{\rm SL}^+$ is related to a graph Laplacian.
So we first introduce some graph-theoretic notation, 
already used in \cite{Kn4} for finite graphs. In the present context
infinite graphs like regular trees arise. See, {\em e.g.} Mohar and 
Woess \cite{MW} for an overview of Laplacians on infinite graphs,
and Fig\`a-Talamanca and Nebbia \cite{FTN} for trees.

For unoriented graphs
$(V,E)$ with bounded vertex degree, we consider the Hilbert space 
$\cH=\cH_V\oplus \cH_{E}$ with  vertex Hilbert space $\cH_V := \ell^2(V)$,
\[{\bf E} := \big\{(v,w) \mid \{v,w\} \in E \big\}
\mbox{ and } \cH_{E}:= \big\{f\in \ell^2({\bf E})\mid f((w,v))=-f((v,w))\big\},\]
with inner product
$\LA f,g\RA_E := \eh\sum_{{\bf e}\in {\bf E}}{f}({\bf e})  \bar g({\bf e})$.

Then the adjoint of 
\[d:\Hi_V\ar \Hi_E,\qquad df\big((v,w)\big):=f(w)-f(v)\]
equals
\[d^*:\Hi_E\ar \Hi_V,\qquad d^*g(v) = -\sum_{w:\,(v,w)\in{\bf E}}g\big((v,w)\big).\]
One defines the operator $\Delta:\Hi\ar\Hi$ by 
\[\Delta := d^*d\oplus dd^* = Q^2 \qmbox{with}Q:=\bsm 0&d^*\\ d&0\esm,\]
so that $\Delta=\Delta_V\oplus \Delta_E$ with {\em vertex Laplacian}
\[\Delta_V f(v)=\sum_{w:\{v,w\}\in E} \big(f(v)-f(w)\big)\qquad (f\in \cH_V)\] 
and {\em edge Laplacian}
\beq
\Delta_E \, g\big((v,w)\big) = 
\sum_{x:(v,x)\in{\bf E}} g\big((v,x)\big) + \sum_{x:(x,w)\in{\bf E}} g\big((x,w)\big)
\qquad (g\in \cH_E).
\Leq{edge:Lap}
Under the assumption of bounded degree $Q$ and thus $\Delta$ 
are bounded, self-adjoint operators, with 
operator norm $\|\Delta\| \le 2\sup_{v\in V} {\rm deg}(v)$. \\
In fact $\big(\Delta,\idty\oplus (-\idty),Q \big)$
is a {\em supersymmetric triple} ({\em i.e.} all operators are self-adjoint with 
$\Delta=Q^2$ and 
$(\idty\oplus (-\idty))^2=\idty$, see, {\em e.g.} Sect.\ 6.3 of 
Cycon, Froese, Kirsch, and Simon \cite{CFKS}).
Apart from zero eigenvalues, the spectra of $\Delta_E$ and $\Delta_V$
coincide (including multiplicities, when finite), since there $Q$ is invertible, with
\[\Delta_E\, d=d\,\Delta_V\qmbox{and}d^*\,\Delta_E=\Delta_V\, d^* .\]
\begin{remarks} \label{rem:I}
\begin{enumerate}
\item 
For the bipartite graphs $(V,E)$ with vertex set $V=V_+\cup V_-$
considered in this article, for $\ell^2(E)$ with counting measure the map 
\[{\cal I}:\ell^2(E)\to\Hi_E\qmbox{,} 
{\cal I}(f)\big((v^\mp,v^\pm)\big):=\pm f\big(\{v^-,v^+\}\big)\quad(v^\pm\in V^\pm)\]
is an isomorphism of Hilbert spaces. We will use ${\cal I}$ implicitly.
\item 
For connected graphs  $(V,E)$, shortest paths define metrics
\[ {\rm dist}\equiv {\rm dist}_V : V\times V \ar \bN_0\qmbox{and} 
{\rm dist}_E: E\times E \ar\bN_0.\]
\end{enumerate}
\end{remarks}
We now determine analytically the spectrum ${\rm spec}(\hT_{\rm SL}^+)$.
This is a proper subset of the set in $\bC$ which 
was determined algebraically in Proposition \ref{prop:algebraic}.
\begin{prop}  \label{prop:SL}
${\rm spec}(\hT_{\rm SL}^+) = \{-\eh,\eh\}\cup C$.
The multiplicities of the  eigenvalues $\pm\eh$ are infinite.
The spectrum of  $\hT_{\rm SL}^+$ on the circle
$C$ is absolutely continuous.
\end{prop}
{\bf Proof.} 
$\bullet$
We first construct the relevant graph $(V,E)$ for $\hT_{\rm SL}^+$.\\
The orbits generated by $-R^{-1}L=\bsm -1&-1\\ 1&0\esm$
are of cardinality three in ${\rm SL}(2,\bZ)$, and similarly for
$-LR^{-1} =\bsm 0&-1\\ 1&-1\esm$. 

So the images of these orbits in ${\rm PSL}(2,\bZ)={\rm SL}(2,\bZ)/\{\pm\idty\}$,
(that is, the orbits $v_+$ of $\pm R^{-1}L$, and the orbits $v_-$ of $\pm LR^{-1}$)
are of order three, too.

The vertex set of the graph is $V:=V_+\cup V_-$, 
with the disjoint sets $V_\pm$ of orbits $v_\pm\subseteq {\rm PSL}(2,\bZ)$.

The edge set $E$ is defined as the set of $\{v_+,v_-\}$
with $v_\pm\in V_\pm$ and $v_+\cap v_-\neq \es$.
For $\{v_+,v_-\}\in E$ the intersection $v_+\cap v_-$
consists of a unique element of the group ${\rm PSL}(2,\bZ)$. 

Taken together, $\pm R^{-1}L$ and $\pm LR^{-1} =\pm J \,R^{-1}L \,J $ 
generate an index two subgroup
of ${\rm PSL}(2,\bZ)$, consisting of the elements which (mod 2)
are in the cyclic subgroup 
\beq
\{\bsm 1&0\\0&1\esm, \bsm 0&1\\1&1\esm, \bsm 1&1\\1&0\esm\}\cong C_3
\qmbox{of} {\rm SL}(2,\bZ/2\bZ).
\Leq{c3}
The coset space disjoint from this subgroup arises by multiplication 
with $\pm J$ and corresponds (mod 2) to the
set $\{\bsm 0&1\\1&0\esm, \bsm 1&1\\0&1\esm, \bsm 1&0\\1&1\esm\}$
of ${\rm SL}(2,\bZ/2\bZ)$--matrices.

Since  ${\rm PSL}(2,\bZ)\cong C_2*C_3$ is a free product, the graph 
$(V,E)$ is thus a disjoint union of two copies of the three-regular tree.\\
$\bullet$
We now use the above identification $E\cong {\rm PSL}(2,\bZ)$.
Dropping the index ${\rm SL}$, on $\cH^+$ the operator 
\beq
\Delta_E := 3(\idty-\hB)=3(\hY_++\hY_-)
=3\idty-\hI+\hR^{-1}\hL + \hL^{-1}\hR + \hL\hR^{-1} + \hR\hL^{-1}
\Leq{Delta:E}
is the edge Laplacian for the graph $(V,E)$
since $(\hR^{-1}\hL)^3 =(\hL^{-1}\hR)^3 =  \hI$
and $\hI^+=\idty_{\cH^+}$ (note that $g\big((v,w)\big)$ appears twice on 
the right hand side of (\ref{edge:Lap})).\\
$\bullet$
The subspace $\cR_{\rm SL}^+$ in (\ref{KRN:SL}) is zero-dimensional, since
it consists of the constant functions that are in $\ell^2({\rm SL}(2,\bZ))$.
So by Lemma \ref{lem:algebraic2}, $\pm1$ are not eigenvalues of 
$\hT_{\rm SL}^+$. We will see that $\pm1\not\in {\rm spec}(\hT_{\rm SL}^+)$, 
when analyzing the operator on $\cN_{\rm SL}^+$.\\
$\bullet$
Thus, using Lemma \ref{lem:algebraic2}, the multiplicities of the eigenvalues 
$\mp \eh$ of $\hT_{\rm SL}^{+,\pm}$ equal the dimensions of 
the $\hY_\pm$--kernel subspaces $\cK_{\rm SL}^{+,\pm}$.

So by (\ref{TJB:formula}) we must construct the eigenfunctions for 
$\hB_{\rm SL}^{+,\pm}$ with eigenvalue~$1$.
This is done by considering the eigenspace for the eigenvalue $0$ of $\Delta_E$. 
For the edge metric ${\rm dist}_E$, and an arbitrary edge 
$g\in E\cong {\rm PSL}(2,\bZ)$, one defines the functions 
$v_{g}\in\ell^2(E)\cong \Hi_E$ (Remark \ref{rem:I}.1),
\beq
v_{g}(h) := \l\{ \begin{array}{ll}
(-2)^{-{\rm dist}_E(g,h)}&\mbox{if $h$ is in the connected component of $g$}\\ 
0& \mbox{otherwise}\end{array}\ri..
\Leq{v:g}
As there are exactly $2^{D+1}$ edges $h\in E$ of distance 
${\rm dist}_E(g,h)=D\in \bN$, the $v_{g}$ are indeed square integrable 
(with $\|v_{g}\|_{E}=\sqrt{3}$).
Notice that $d^* v_{g}=0\in \cH_V$. 

So $v_{g}$ is an eigenfunction with eigenvalue zero of $\Delta_E=d\,d^*$ 
(a harmonic function) which 
does not correspond to an eigenfunction of $\Delta_V$. 

By (anti--)symmetrizing
 $v_g$ w.r.t.\ the $\hJ$--action (that is, w.r.t.\ the two connected components of the
 graph) we show
$\pm\eh\in {\rm spec}(\hT_{\rm SL}^+)$. The multiplicity of these
eigenvalues is infinite, although the 
set $\{v_g\mid g\in {\rm PSL}(2,\bZ)\}$ is linear dependent.\\
$\bullet$
The rest of ${\rm spec}(\hT_{\rm SL}^{+})$ relates to the 
subspace $\cN_{\rm SL}^{+}$. Since $\cR_{\rm SL}^+=\{0\}$, 
this corresponds to ${\rm ran}(\Delta_E)$, 
or by supersymmetry, to ${\rm ran}(\Delta_V)$.
 
It is known (see, {\em e.g.} \cite{FTN}, Chapter II.6) that the spectrum
of the vertex Laplacian on the $(q+1)$--regular tree equals 
\beq
{\rm spec}(\Delta_V) = \big[q+1-2\sqrt{q},\,q+1+ 2\sqrt{q}\,\big] 
\Leq{spec:DV} 
and is absolutely continuous.

However, we repeat the proof for the Laplacian at hand, 
since we need it in the real and the adelic cases below.
Instead of (\ref{spec:DV}) we show the estimate 
\beq
{\rm spec}({\rm Ad}) = \big[-\sqrt{8}, \sqrt{8}\,\big] 
\Leq{spec:Ad} 
for the {\em adjacency matrix} ${\rm Ad}\in \hB(\cH_V)$ which equals
${\rm Ad}\,f(v)=\sum_{w:\{v,w\}\in E}  f(w)$. (\ref{spec:Ad}) is equivalent to 
(\ref{spec:DV}), since $(V,E)$ is three-regular.
We prove that the spectrum is contained in (\ref{spec:Ad}) by constructing 
the resolvent $({\rm Ad}-\lambda \idty)^{-1}$ for all
$\lambda\in\bC\setminus [-\sqrt{8}, \sqrt{8}\,\big]$.
We define the family of operators
\beq
D(k)\equiv D_{{\rm SL}}(k)\in \hB(\cH_V)\qmbox{,}
\big(D(k)f\big)(v):= \sum_{w:\;{\rm dist}(v,w)=k} f(w)\qquad (k\in\bN_0).
\Leq{Def:Dk}
So in particular $D(0)= \idty_{\cH_V}$ and $D(1)={\rm Ad}$. The $D(k)$ 
are self-adjoint.\\
$\bullet$
Their norm is estimated as follows for $k\in\bN$. Given an arbitrary tree root $u\in V$,
one sets $d(v):= {\rm dist}(u,v)$ and orthogonally decomposes the Hilbert space as 
\[\cH_V=\bigoplus_{\ell\in\bZ/k\bZ} \cH_\ell \; \qmbox{with} 
\cH_\ell:= \big\{f\in \cH_V\mid 
\supp f \subseteq C_\ell \big\}\] 
for $C_\ell:= \{v\in V\mid [d(v)]_k = \ell\  \}$, with $[x]_k:=x\ ({\rm mod}\, k)$.\\
Accordingly, $D(k)=(D(k)_{\ell,m})$ with restricted operators 
$D(k)_{\ell,m}: \cH_\ell\ar\cH_m$. Their operator norms are bounded by 
\beq
\|D(k)_{\ell,m}\|\le 3\;2^{k/2}\qquad\big(k\in\bN,\; \ell,m\in\bZ/k\bZ\big).
\Leq{Dk:norm}
To see this, we notice that for vertices $v\in C_\ell$ and $w\in C_m$ with 
${\rm dist}(v,w)=k$ there is a unique vertex $x\in V$ in the 
intersection of the connecting paths
$[u,v]$, $[u,w]$ and $[v,w]\subseteq V$, and 
$2\,d(x)= \ell+m \ ({\rm mod }\, k)$.
By the triangle inequality ${\rm dist}(x,w)\le k$, so given 
$\ell$ and $m$ there occur at most three values of ${\rm dist}(x,w)=d(w)-d(x)$. 
Let $W\equiv W(\ell,m)\subseteq \bN_0$ 
be the set of these values.

The image of $f\in\cH_\ell$ has the squared norm  $\|D(k)_{\ell,m}f\|^2=$
\beqno
&=& 
\sum_{w\in C_m}\Big| \sum_{\stackrel{v\in C_\ell:}{ {\rm dist}(v,w)=k}}f(v) \Big|^2
=\sum_{\widetilde w\in W} 
\sum_{x\in C_{m-[\widetilde w]_k}}
\sum_{\stackrel{w\in C_m:} { {\rm dist}(x,w)=\widetilde w}}
\Big|\sum_{\stackrel{v\in C_\ell:\,{\rm dist}(x,v)
=k-\widetilde w,}{{\rm dist}(v,w)=k}}f(v) \Big|^2\\
&\le&\sum_{\widetilde w\in W} 
\dh\; 2^{\widetilde w}\sum_{x\in C_{m-[\widetilde w]_k}}
\Big|\sum_{\stackrel{v\in C_\ell:\,d(v)\ge d(x), }{\,{\rm dist}(x,v)
=k-\widetilde w}}f(v) \Big|^2\\
&\le& \sum_{\widetilde w\in W} 
\dh\; 2^{\widetilde w}\sum_{x\in C_{m-[\widetilde w]_k}}\dh\; 2^{k-\widetilde w}
\sum_{\stackrel{v\in C_\ell:\,d(v)\ge d(x),}{{\rm dist}(x,v)
=k-\widetilde w}}\big|f(v) \big|^2
\ \le \ \textstyle\frac{27}{4}\, 2^{k}\|f\|^2,
\eeqno
the second to last inequality being H\"older's. This proves (\ref{Dk:norm}).
Since the $L^2$--matrix norm is bounded above by the Hilbert-Schmidt
norm, we thus get 
\[\|D_{{\rm SL}}(k)\|\le \Big(\sum_{\ell,m \in\bZ/k\bZ}\|D(k)_{\ell,m}\|^2 \Big)^{1/2}
\le 3\,k\, 2^{k/2} \qquad(k\in\bN).\]
$\bullet$
For parameters $c,d\in\bC$ the series 
\beq
\cD(c)\equiv \cD_{{\rm SL}}(c) := d\,\idty + \sum_{k=1}^\infty D(k)\, 2^{-ck} 
\Leq{resolvent:SL}
thus converges to a bounded operator, if $\Re(c)>\eh$. \\
Conversely, the estimate for the partial sums
$\|(d\,\idty +\sum_{k=1}^K D(k)\, 2^{-ck})\delta_u\|=$
\beq
=\Big(d^2+\dh   \textstyle\sum_{k=1}^K 2^{(1-2c)k}\Big)^{1/2}
=\Big( \textstyle d^2+3 \frac{2^{(1-2 c) K}\,-1}{2-4^c} \Big)^{1/2}
\Leq{partial:sums} 
for all $K\in\bN$ shows
divergence for $\Re(c)<\eh$, and similarly for $\Re(c)=\eh$.\\
$\bullet$
For $\Re(c)>\eh$ the operator $\cD_{{\rm SL}}(c)$ is a resolvent for ${\rm Ad}$:
\beqno
({\rm Ad}-\lambda \idty)\;\cD(c)
&=& (3\; 2^{-c}-\lambda d)\idty  -\sum_{k=1}^\infty\lambda D(k)2^{-ck} \\
&&+\sum_{k=2}^\infty2D(k-1)\,2^{-ck} 
+\sum_{k=0}^\infty D(k+1)\, 2^{-ck} \\
&=&(3\; 2^{-c}-\lambda d)\idty +
\sum_{k=1}^\infty  \big(2^{1-c}-\lambda +2^{c} \big)D(k)2^{-ck} =\idty
\eeqno
for $d:=\frac{3\,2^{-c}-1}{\lambda}$ and $2^{1-c} +2^{c} 
=\lambda$ or $c_\pm(\lambda):=\log_2 \big(\frac{1}{2} (\la\pm\sqrt{\la^2-8})\big)$.
Since $\max\big(\Re(c_\pm(\lambda))\big) > \eh$ if and only if 
$\lambda\in\bC\setminus [-\sqrt{8}, \sqrt{8}\,\big]$, 
this proves (\ref{spec:Ad}), and thus (\ref{spec:DV}) for $q=2$.\\
$\bullet$
In our case the spectral interval 
\[{\rm spec}(\Delta_E)\setminus \{0\}
= {\rm spec}(\Delta_V) = \big[ 3 - \sqrt{8}, 3 + \sqrt{8} \big] \] 
is, using (\ref{TJB:formula}) and $\Delta_E=3(\idty_{\cH^+}-\hB^+)$, 
the image of the circle $C$ under the maps 
\[x\;\longmapsto\; 3\pm (2x+1/x).\] 
The sign depends on the  $\hJ^+$ subspace.
Each point $y\in(3-\sqrt{8},3+ \sqrt{8})$
has the two complex conjugate preimages 
\beq
 \ev \left( 3-y \pm \imath\sqrt{6y-y^2-1}\right)\mbox{ respectively }
 -\ev \left( 3-y \pm \imath\sqrt{6y-y^2-1}\right).
 \Leq{double:cov}
As for all representations considered in this article, 
the spectrum of $\hT_{\rm SL}^{+}$ is invariant 
under complex conjugation. So $C$ belongs to ${\rm spec}(\hT_{\rm SL}^{+})$,
but the only real points of the spectrum are $\pm\eh$ and $\pm 1/\sqrt{2}$.\\
$\bullet$
As the inverse transformations (\ref{double:cov}) are analytic
and non--constant on the open interval $(3-\sqrt{8},3+ \sqrt{8})$, 
absolute continuity of ${\rm spec}(\Delta_V)$ leads to
absolute continuity of ${\rm spec}(\hT_{\rm SL}^{+})$.
\hfill $\Box$\\[2mm]

We are particularly  interested in the contractions $\hT_S^+$, for the 
real 'infinite' case $S=\{\infty\}$, the finite adele case $S= \bP$ 
and the adelic case $S= \bP_\infty$,
which we now analyze in succession.

The first steps will always consist in considering the orbits of 
${\rm SL}(2,\bZ)$ in $\bZ_S^2$,
since ${\rm SL}(2,\bZ)$ is generated by left and right addition 
$L=\bsm 1&1\\0&1\esm$, $R=\bsm 1&0\\1&1\esm$, and their inverses.
%
\subsection{The Real Case}\label{subs:rc}
%
We begin with the simplest of our three sets of places, 
$S=\{\infty\}$ and thus analyze $\hT_\infty$ on $\ell^2(\bZ^2)$. 
The ${\rm SL}(2,\bZ)$--orbits in $\bZ^2$ are of the form 
\beq
n\,\Lambda\quad\mbox{with }\Lambda:=\{\bsm a\\ b\esm\in \bZ^2\mid \gcd(a,b)=1\}
\mbox{ and }n\in\bN_0.
\Leq{SL:orbits:real}
Except for $0\,\Lambda=\{\bsm0\\0\esm\}$, the actions on these orbits 
are mutually isomorphic. Thus in the orthogonal decomposition 
\beq
\ell^2\big(\bZ^2\big)=
\bigoplus_{n\in\bN_0}\ell^2(n\,\Lambda)\mbox{ , and }
\hT_\infty= \bigoplus_{n\in\bN_0}\hT_{\infty,n} \mbox{ with }
\hT_{\infty,n}:=\hT_{\infty}\rstr_{\ell^2(n\,\Lambda)}
\Leq{split:real}
we need only consider the trivial case $n=0$, and $n=1$, setting 
\[\cH_\Lambda := \ell^2(\Lambda)\qmbox{and}\hT_\Lambda := \hT_{\infty,1}.\]
The action of ${\rm SL}(2,\bZ)$ on $\Lambda$ is isomorphic to 
its left action on the coset space by the parabolic subgroup
$L^\bZ=\{\bsm 1&n\\0&1\esm\mid n\in\bZ\}$, via
\beq
{\rm SL}(2,\bZ)/ L^\bZ\ \cong \ \Lambda
\qmbox{,}\l[\bsm a&b\\c&d\esm\ri]\mapsto \bsm a\\ c\esm\ ,
\Leq{left:reg}
fitting with the left regular representation.

Thus the action of ${\rm SL}(2,\bZ)$ on $\Lambda$ is not free, 
unlike the one on itself.
However we will see in Lemma \ref{lem:F} that in a precise sense 
it is not far from being free, 
a fact that is vital for the spectrum of~$\hT_\Lambda$.

The typical fiber of the bundle 
\beq
\Pi:{\rm SL}(2,\bZ)\ar\Lambda
\Leq{def:Pi}
induced by (\ref{left:reg}) is isomorphic to $\bZ$.
We define a section $\Lambda\ar{\rm SL}(2,\bZ)$ 
by choosing for $\bsm a\\ c\esm\in\Lambda$ the element 
$\bsm a&b-ka\\c&d-kc\esm$ of the fiber over $\bsm a\\ c\esm$
with minimal $k\in\bZ$ so that 
$\LA \bsm a\\ c\esm,\bsm b-ka\\d-kc\esm\RA\le 0$.
Thus we obtain a trivialization 
$\Pi\times F:{\rm SL}(2,\bZ)\ar\Lambda\times \bZ$
of the bundle with second factor
\beq
F:{\rm SL}(2,\bZ)\ar \bZ \qmbox{,} 
\bsm a&b\\c&d\esm\mapsto \lceil(ab+cd)/(a^2+c^2)\rceil,
\Leq{def:F}
(the argument of the ceil function being the unipotent parameter of the 
Iwasawa decomposition of $\bsm a&b\\c&d\esm\in {\rm SL}(2,\bR)$). 

The matrices appearing in the following lemma are
$\pm R^{-1}L$, $\pm L^{-1}R$, $\pm RL^{-1}$ and $\pm LR^{-1}$. 
Apart from the identity, it is their unitary representations that are the terms
of the operator $\hB^+$, see \eqref{Delta:E}.
\begin{lemma}\label{lem:F}
For $A\in {\rm SL}(2,\bZ)$ 
\[ | F(MA)-F(A) | \le 1 \qmbox{if} 
M\in\pm \big\{
\bsm  -1 &-1\\ 1&0\esm, 
\bsm  0 &1\\ -1&-1\esm,
\bsm  -1 &1\\ -1&0\esm,
\bsm  0 &-1\\ 1&-1\esm\big\}. \] 
\end{lemma}
{\bf Proof.}
We set $A:=\bsm a&b\\c&d\esm$ and $\tilde F(A) := (ab+cd)/(a^2+c^2)$. 
Since $\tilde F$ is even, we need only consider the positive sign in the list of $M$.
These $M$ are elliptic of order three, with 
$\bsm  0 &1\\ -1&-1\esm = \bsm  -1 &-1\\ 1&0\esm^{-1}$
and $\bsm  0 &-1\\ 1&-1\esm  = \bsm  -1 &1\\ -1&0\esm^{-1}$.\\
$ \big| \tilde F(\bsm  -1 &-1\\ 1&0\esm A)-\tilde F(A)  \big|
=|(a^2-a c-c^2)/((a^2+c^2)(2a^2+2a c +c^2))|$ is smaller than 
$\frac{\sqrt{5}}{2}/(a^2+c^2)$ and thus $<1$ if $a^2+c^2>1$. 
It is $\le 1$ if $a^2+c^2=1$.\\
The case $\big| \tilde F(\bsm  -1 &1\\ -1&0\esm A)-\tilde F(A)  \big|
= |(a^2+a c-c^2)/((a^2+c^2)(2a^2-2a c +c^2))|$ is similar.
\hfill $\Box$\\[2mm]
According to  (\ref{splitting}) and (\ref{KRN}) we split the 
Hilbert space $\cH_\Lambda$ into the orthogonal  subspaces 
\beq
\cH_\Lambda=\cH_\Lambda^+\oplus\cH_\Lambda^-\qmbox{with}
\cH_\Lambda^+=\cK_\Lambda^+\oplus\cR_\Lambda^+\oplus\cN_\Lambda^+.
\Leq{KRN:R}
\begin{prop} \label{prop:real} %
${\rm spec}\big(\hT_{\infty,0}\big) = \{1\}$, 
${\rm spec}\big(\hT_{\infty,n}\big) = {\rm spec}\big(\hT_\Lambda\big)$ $(n\in \bN)$, 
and
\[{\rm spec}\big( \hT_\Lambda^+ \big)= \{-\eh,\eh\}\cup C .\]
Like for the case $\hT_{\rm SL}^+$, the spectrum of  $\hT_\Lambda^+$ on 
$C$ is absolutely continuous. 
\end{prop}
{\bf Proof.} 
$\bullet$ Since $\ell^2(0\,\Lambda)\cong\bC$, the first statement is obvious.\\
$\bullet$
The pull-back $\ell^2(n\,\Lambda)\ar \ell^2(\Lambda)$ induced by multiplication
$\Lambda \ar n\,\Lambda$ by $n$ is unitary for $n\in\bN$
and then conjugates $\hT_{\infty,n}$ with $\hT_\Lambda$.\\
$\bullet$
We first fix the function spaces involved.
The operator induced by the map $\Pi$ from (\ref{def:Pi}) is 
defined in two steps. First we define on the subspace 
\[\widetilde U:= \big\{\phi\in \bC^{{\rm SL}(2,\bZ)}\mid
\forall x\in\Lambda:\phi\rstr_{\Pi^{-1}(x)} \in\ell^1\big(\Pi^{-1}(x)\big)\big\}\]
with absolutely integrable fibers the operator
\beq\widetilde \Pi:\widetilde U\ar  \bC^{\Lambda}\qmbox{,}
\big(\widetilde \Pi \phi\big)(x)= \sum_{y\in  \Pi^{-1}(x)} \phi(y).
\Leq{def:whPi}
Then we set 
\[U:= \big\{\phi\in\widetilde U\mid \widetilde \Pi(|\phi |)\in \ell^2(\Lambda)\big\}
\qmbox{and}\widehat \Pi:= \widetilde \Pi\rstr_U.\] 
$U$ is a subspace of $\ell^2\big({\rm SL}(2,\bZ)\big)$, since 
$\sum_{y\in \widehat \Pi^{-1}(x)} |\phi(y)|^2\le 
(\sum_{y\in \widehat \Pi^{-1}(x)} |\phi(y)|)^2$.

$\widehat \Pi$ defines an unbounded operator between Hilbert spaces,
since the fibers of the bundle projection
$\Pi$ are isomorphic to $\bZ$, {\em i.e.} infinite.
Nevertheless, on certain subspaces this norm is finite: 
For all $N\in\bN_0$, with $F$ defined in (\ref{def:F})
\[U_N:=\big\{\phi\in \ell^2\big({\rm SL}(2,\bZ)\big)\mid 
F(\supp\;\phi)\subseteq \{-N,\ldots,N\}\big\}\ \subseteq \ U,\]
and by H\"older's inequality
\beq
\big\|\widehat \Pi\phi\big\|\leq \sqrt{2N+1}\,\|\phi\|\qquad(\phi\in U_N).
\Leq{UN:estimate}
So $U$ is a dense subspace of $\ell^2\big({\rm SL}(2,\bZ)\big)$.\\
$\bullet$
Conversely, the embedding 
$E:\Lambda\cong \Lambda\times \{0\}\hookrightarrow {\rm SL}(2,\bZ)$ defined by 
the bijection $\Pi\times F:{\rm SL}(2,\bZ)\ar\Lambda\times \bZ$ 
induces an isometric embedding 
\[\widehat E:\ell^2(\Lambda)\ar U
\qmbox{,} 
\big(\widehat E\phi\big)(y)
=\l\{\begin{array}{ll}\phi\big(E^{-1}(y)\big)&\mbox{if }F(y)=0\\
0&\mbox{else}\end{array}\ri.\]
of Hilbert spaces (and $\widehat E(\ell^2(\Lambda))= U_0$).\\
$\bullet$
We begin with the subspace $\cK_\Lambda^+$ in (\ref{KRN:R}).
Due to Lemma \ref{lem:algebraic2}
the spectra of $\hT_{\Lambda}^+$ on the $J$--subspaces $\cK_\Lambda^{+,\pm}$
equal $\{-\eh\}$ respectively $\{\eh\}$.
$\cK_\Lambda^{+,\pm}$ are non-trivial, since the eigenfunctions 
$v_g$ of $\hT^+_{\rm SL}$, defined in (\ref{v:g})
are in the domain $U$, but not in the kernel of the projection operator $\widehat \Pi$:
\begin{enumerate}[-] 
\item
It is absolutely summable over the fibers of (\ref{def:Pi}) at any 
$h\in {\rm SL}(2,\bZ)$, since 
$\sum _{m=-\infty }^{\infty } 2^{-|m|} = 3< \infty$ for $h=g$, and similar else.
\item
It is non-vanishing, since the projection at $g$ is 
$\big(\widehat \Pi v_g\big)(g)
= \mbox{\footnotesize ${\displaystyle\sum _{m\in\bZ}}$} (-2)^{-|m|} = \ed$.
\end{enumerate}
$\bullet$
To show that the spectrum of $\hT^+_\Lambda$ on $\cN^+_\Lambda$ 
is contained in the
circle $C$, we use (\ref{TJB:formula}) and prove that 
${\rm spec}(\hB^+_\Lambda\rstr_{\cN^+_\Lambda})\subseteq 
\big[-\sqrt{8}/3,\,\sqrt{8}/3\,\big]$.
For that we compare $\hB^+_\Lambda$ with $\hB^+_{\rm SL}$, whose resolvent on
$\cN^+_{\rm SL}$ was given in (\ref{resolvent:SL}). 
$\hB_{\Lambda}^+$ on the Hilbert space
$\cN^+_\Lambda$ is analyzed with the help
of the commutative diagram (indeed, $\hB_{\rm SL}$ restricts to $U$)
\beq
\begin{CD}
U @>\hB_{\rm SL}>>U\\
@A\widehat E AA@VV\widehat  \Pi V\\
\ell^2(\Lambda) @>\hB_\Lambda>> \ell^2(\Lambda)
\end{CD}\ .
\Leq{CD1}
The operators $D_{{\rm SL}}(k)$ from (\ref{Def:Dk}) map 
by Lemma \ref{lem:F} the subspace
$U_0$ into $U_k$.
Therefore, by (\ref{UN:estimate}) for $\Re(c)>\eh$ the resolvent 
$\cD_{{\rm SL}}(c)=\big(\hB^+_{\rm SL}-\lambda(c)\idty\big)^{-1}$ of 
$\hB^+_{\rm SL}$ has the property
\[\cD_{{\rm SL}}(c)(U^+_0)\subseteq U^+\qmbox{, and}
\big(\hB^+_\Lambda-\lambda(c)\idty\big)^{-1} =
\widehat  \Pi \circ \big(\hB^+_{\rm SL}-\lambda(c)\idty\big)^{-1}\circ\widehat E \]
is bounded. This shows that 
${\rm spec}\big( \hB_\Lambda^+\rstr_{\cN^+_\Lambda} \big)
\subseteq [-\sqrt{8}/3,\sqrt{8}/3]$.
Prop.\ \ref{prop:ABT}  then implies that
${\rm spec}\big( \hT_\Lambda^+\rstr_{\cN^+_\Lambda} \big)\subseteq C$.\\
$\bullet$
The converse inclusion 
${\rm spec}\big( \hT_\Lambda^+\rstr_{\cN^+_\Lambda} \big)\supseteq C$
is provided by (\ref{partial:sums}), together 
with Lemma \ref{lem:F}. They imply divergence of the partial sums
for the resolvent $\big(\hB_\Lambda^+ - \lambda(c)\idty\big)^{-1}$ in the case  
$\Re(c) < \eh$.  \\
$\bullet$
The absolute continuity of 
${\rm spec}\big( \hB_\Lambda^+\rstr_{\cN^+_\Lambda} \big)$
follows from absolute continuity of 
${\rm spec}\big( \hB_{\rm SL}^+\rstr_{\cN^+_{\rm SL}} \big)$; the one of
${\rm spec}\big( \hT_\Lambda^+\rstr_{\cN^+_\Lambda} \big)$ then follows using
the general formula~(\ref{rep:T}).
\hfill$\Box$
%

\subsection{Representations of ${\rm SL}(2,\bZ/n\bZ)$}\label{subs:ZnZ}

We denote the binary modular congruence groups by
\[{\rm G}_n:={\rm SL}(2,\bZ/n\bZ)\qquad(n\in\bN).\] 
Via the surjective homomorphism
${\rm SL}(2,\bZ)\ar {\rm G}_n$ the group ${\rm SL}(2,\bZ)$ acts on 
the Hilbert space
\[h(n):=\ell^2\big((\bZ/n\bZ)^2\big)\]
by the permutation representation. As remarked above, the ${\rm G}_n$--orbit
\beq
\Lambda(n)=\big\{\bsm a\\b\esm \in(\bZ/n\bZ)^2\mid 
\gcd(a,b,n)=1\big\}\qquad(n\in\bN)
\Leq{def:Lambda:n}
has cardinality 
$|\Lambda(n)| = J_2(n) = n^2\prod_{p\in\bP:\, p|n}(1-p^{-2})$. 
To it we associate the $|\Lambda(n)|$--dimensional, 
${\rm SL}(2,\bZ)$-invariant  subspace
\[\tilde h(n):=
\big\{f\in h(n)\mid \supp(f)\subseteq \Lambda(n)\}.\]
We use small bold letters to abbreviate operators on $h(n)$. 
So with $\hT$ from ({\ref{the:op})  
\beq
\htt(n):= \hT_{(\bZ/n\bZ)^2}\in\hB\big(h(n)\big)
\qmbox{and}\tilde t(n):= \hT_{\Lambda(n)}\in\hB\big(\tilde h(n)\big).
\Leq{tn:tildetn}
The spectral theory of the operators $\tilde \htt(n)$ on $\tilde h(n)$  
is related to the representation theory of the group ${\rm G}_n$,
since the action of ${\rm SL}(2,\bZ)$ on $\Lambda(n)$ gives rise to 
a representation of ${\rm G}_n$ on $\tilde h(n)$. 

Unitary Fourier transform is denoted by
\[\cF_n: h(n)\to h(n)\mbox{ , }\ 
(\cF_nf)(\ell)=n^{-1}\!\!\!\sum_{k\in (\bZ/n\bZ)^2}\!\!\!
\exp(-2\pi\imath \LA k,\ell\RA/n)f(k)\ \quad (n\in\bN).\]
\begin{lemma}  \label{lem:commu}
On the Hilbert spaces $h(n)$, with ${\bf j}(n) := \hJ_{(\bZ/n\bZ)^2}$
\[ [\cF_n,{\bf j}(n)]=0\qmbox{,} 
[\cF_n,{\bf b}(n)]=0\qmbox{and}[{\bf j}(n)\cF_n,\htt(n)]=0.\]
\end{lemma}
\textbf{Proof.}
$\bullet$
In general for $O\in{\rm SL}(2,\bZ)$ represented by ${\bf o}(n)$ we get 
\beqno
\cF_n {\bf o}(n) f(k)&=&n^{-1}\sum_{\ell\in (\bZ/n\bZ)^2} 
\exp\big(-2\pi\imath \LA k,\ell\RA/n\big) f\circ O^{-1}(\ell)\\
&=&n^{-1}\sum_{m\in (\bZ/n\bZ)^2} 
\exp\big(-2\pi\imath \LA k,O(m)\RA/n\big) f(m)\\
&=&n^{-1}\sum_{m\in (\bZ/n\bZ)^2} 
\exp\big(-2\pi\imath \LA O^\top(k),m\RA/n\big) f(m)\\
&=&
(\cF_n f)\big(O^\top(k)\big) = 
\big(({\bf o}(n)^{\mbox{\scriptsize\boldmath $\top$}})^{-1}\cF_nf\big)(k).
\eeqno
$\bullet$
As $(J^\top)^{-1}=J$, $[\cF_n,{\bf j}(n)] = 0$, \\
$\bullet$
and as $L^\top =R$, $\cF_n$ permutes the representations 
${\bf l}(n)$ and ${\bf r}(n)^*$.
So $\cF_n$ conjugates $\htt(n)$ and $\htt(n)^*$. 
But the same does ${\bf j}(n)$, see (\ref{identi}).
So $[{\bf j}(n)\cF_n,\htt(n)] =0$. \\
$\bullet$
The case of ${\bf b}(n)\equiv\hB=
\ed(\hI-\hR^{-1}\hL - \hL^{-1}\hR - \hL\hR^{-1} - \hR\hL^{-1})$ is similar.
\hfill $\Box$\\[2mm]
${\bf j}(n)\cF_n$ is not a multiple of the identity if $n\in\bN\setminus\{1\}$,
although  both ${\bf j}(n)$ and $\cF_n$ have the same square ${\bf i}(n)$.
For an eigenfunction $f\in h(n)$ of $\htt(n)$ by Lemma \ref{lem:commu} 
we either obtain a linearly independent eigenfunction 
${\bf j}(n)\cF_nf$ or a non-trivial symmetry of~$f$. 
The second alternative applies, {\em e.g.}, to 
the mean zero eigenfunctions $f\in\tilde h(n)$ of $\tilde t(n)$ if $n\in \bP$.

For the prime case $n\in \bP$ the spectral theory of the operators $\tilde t(n)$
is partly done in \cite{Kn4}.
So we review here the theory of ${\rm G}_n$ representations  for 
arbitrary $n\in\bN$ and then apply it to the operators~$\tilde t(n)$.

The following example shows that the operators $\htt(n)$ need not be semisimple.
\begin{example} \label{ex:not:semisimple}
The operator $\htt(6)$ is defective.
The vector $v\in h^-(6)$, 
\[v:=\bsm 
 0  & 1  & 0 & 0  & 0 & -1 & \
 -1 & 0  & 0 & 0  & 0 & 0  & \
 0  & 0  & 0 & -1 & 0 & 1  & \
 0  & -1 & 0 & 0  & 0 & 1  & \
 0  & -1 & 0 & 1  & 0 & 0  & \
 1  & 0  & 0 & 0  & 0 & 0
 \esm\] 
(in lexical order of $(\bZ/6\bZ)\times(\bZ/6\bZ)$ in the least residue system 
modulo 6), is a generalized $\htt(6)$-eigenvector of eigenvalue 0, withÊ
$\htt(6)^2v=0$, but  $\htt(6) v\neq0$. \hfill $\Diamond$

\end {example}

\subsubsection{The Graph of the Regular ${\rm G}_n$--Representation}
By general wisdom the regular representation of ${\rm G}_n$
contains all irreducible representations, with multiplicity
given by cardinality of conjugacy classes. The operator $\hT_{{\rm G}_n}$
on $\ell^2({\rm G}_n)$ thus restricts to these subrepresentations.

For $n>2$ the operator $\hT^+_{{\rm G}_n}$ is related to the Laplacians
of a graph $(V_n,E_n)$ in a way that is analogous to
the case of  $\hT^+_{{\rm SL}}$ treated in Section \ref{subs:rr}.
\begin{itemize}
\item 
Again the bipartite vertex set $V_n=V_n^+\cup V_n^-$ is composed of the 
set $V_n^+$ of orbits of $-R^{-1}L\in{\rm G}_n$ and
the orbit set $V_n^+$ of  $-L^{-1}R\in{\rm G}_n$.
\item 
The edge set $E_n$ equals the group, the 
edges ${v^-,v^+}$ connecting orbits $v^\pm$ with $v^-\cap v^+=\{g\}$
with $g\in {\rm G}_n$. 
\end{itemize}
Figure \ref{fig:borel} shows the graph of the group ${\rm G}_5$.
\begin{prop}  %
\ For $n\in\bN$, $n\ge3$ with prime decomposition 
$n=\prod_{i=1}^s p_i^{k_i}$ \\
$\bullet$
the graph $(V_n,E_n)$ is three-regular and thus 
has $|V_n|=\frac{2}{3}|E_n|$ vertices, with the number 
$|E_n|=|{\rm G}_n|=\prod_{i=1}^s p_i^{3k_i-2}(p_i^2-1)$ of edges.\\
$\bullet$
For odd $n$ it is connected. If $n/2$ is odd, it has two components and
for $4|n$ it has four components. All components are isomorphic.\\
$\bullet$
The girth of a graph $(V,E)$ being the length of a shortest cycle,
\[{\rm girth}\big((V_n,E_n)\big)\ge
2 \left\lfloor \frac{\cosh ^{-1}\left(\frac{\sqrt{5}
   n}{4}\right)}{{\rm csch}^{-1}(2)}\right\rfloor -2
\sim 2\log_\Phi(n)\qquad\mbox{with golden ratio }\Phi.\]
\end{prop}
\textbf{Proof.}
$\bullet$
For $n\ge 3$ the matrices 
$-R^{-1}L =\bsm-1&-1\\1&0\esm$, $(-R^{-1}L)^2 = \bsm 0&1\\-1&-1\esm$, 
$-LR^{-1} =\bsm 0&-1\\1&-1\esm$ and $(-LR^{-1})^2 =\bsm -1&1\\-1&0\esm$
are all different in ${\rm G}_n$.
Thus all orbits are of order three, and $(V_n,E_n)$ is three-regular.\\
$\bullet$
The group ${\rm GL}(2,\bZ/n\bZ)$ of invertible matrices over the
residue class ring $\bZ/n\bZ$ is
isomorphic to  
${\rm GL}(2,\bZ/p_1^{k_1}\bZ)\times \ldots\times {\rm GL}(2,\bZ/p_s^{k_s}\bZ)$.
Moreover,  
\[|{\rm GL}(2,\bZ/p^k\bZ)|=p^{4(k-1)}|{\rm GL}(2,\bZ/p\bZ)|,
\mbox{ with }  |{\rm GL}(2,\bZ/p\bZ)|=(p^2-1)(p^2-p),\]
the formula going back to Camille Jordan. 
${\rm SL}(2,\bZ/n\bZ)$
is the kernel of the determinant homomorphism of ${\rm GL}(2,\bZ/n\bZ)$.
So with $\varphi(p^k)=(p-1)p^{k-1}$
the cardinality equals
\beq
|{\rm SL}(2,\bZ/p^k\bZ)| = |{\rm GL}(2,\bZ/p^k\bZ)|/\varphi(p^k) = p^{3k-2}(p^2-1).
\Leq{dim:SL}
We also have the group isomorphism
\beq
{\rm SL}(2,\bZ/n\bZ)\cong 
{\rm SL}(2,\bZ/p_1^{k_1}\bZ)\times \ldots\times {\rm SL}(2,\bZ/p_s^{k_s}\bZ),
\Leq{congo} 
since the 
above homomorphism for ${\rm GL}(2,\bZ/n\bZ)$ respects the determinant $1$.\\
$\bullet$
We first show that for $2|n$ there are at least two graph components in 
$(V_n,E_n)$.
Then the homomorphism ${\rm G}_n\to {\rm G}_2$ shows, like in (\ref{c3}),
that the identity $\bsm 1&0\\0&1\esm\in {\rm G}_n$ lies in a component projecting
onto the cyclic subgroup $C_3$ of ${\rm G}_2$ of index 2, missing 
the coset space of $J=\bsm 0&1\\-1&0\esm$.

For $4|n$ the homomorphism ${\rm G}_n\to {\rm G}_4$ shows additionally
that the  image of the identity component does not contain 
$I=\bsm -1&0\\0&-1\esm=J^2\in {\rm G}_4$. Unlike in ${\rm G}_2$, this
is not the identity. So the graph has at least  four components, corresponding
to the powers of $J$.

The graph $(V,E)$ constructed in the proof of Proposition \ref{prop:SL}
has the group ${\rm PSL}(2,\bZ)$ as edge set $E$, is the disjoint union of 
two copies of the three-regular tree, and thus is covered by
a graph with edge set ${\rm SL}(2,\bZ)$, consisting of four three-regular trees.
As the latter graph covers $(V_n,E_n)$, the latter has no more components.\\ 
$\bullet$
Multiplication by $J^k$ provides graph isomorphisms of the components.\\
$\bullet$
The norms of the generators $-R^{-1}L$, $-LR^{-1}\in {\rm SL}(2,\bZ)$ 
and their inverses equal the golden ratio $\Phi=\eh (\sqrt{5}+1)$. 
So all their $k$-fold products $X=\bsm x_{1,1}&x_{1,2}\\ x_{2,1}&x_{2,2}\esm$ 
have norms $\|X\|\le \Phi^k$.
By induction $\max(|x_{i,j}|)\le f_{k+1}$, $f_\ell$ being the $\ell$th
Fibonacci number, and the next to maximal $|x_{i,j}|\le f_k$. So for 
values of $k$ with $f_k\le n/2$, all $X$ are different ${\rm mod}\ n$. 
Solving for $k$, using Binet's formula $f_k=(\Phi^k-(-\Phi)^{-k})/\sqrt{5}$ 
gives the result.
\hfill $\Box$

\subsubsection{Representations of 
${\rm SL}(2,\bZ/n\bZ)$ for Primes {\large $n$}}
We begin with the case $n=p\in\bP$ of primes, with the
field $\bF_p=\bZ/p\bZ$. 
It is then known (see Naimark and \v Stern \cite{NS},  Section II, \S5, and 
Lafferty and Rockmore \cite{LR}, Section 2) that the 
irreducible representations of 
${\rm SL}(2,\bF_p)$ are divided into the classes of principal respectively discrete
(or cuspidal) series. Denoting by 
\[{\rm B}_n:= \big\{\!\bsm a&b\\0&a^{-1}\esm\in {\rm G}_n \big\} \qmbox{and}
{\rm U}_n := \{\bsm1&b\\0&1\esm\in {\rm B}_n \} \qquad(n\in \bN)\]
the Borel respectively unipotent subgroups, the discrete series representations 
of ${\rm G}_p$ are characterized by the property that their restriction to 
${\rm U}_p$ does not contain the trivial representation.

The principal series representations of ${\rm G}_p$ 
are the irreducible subrepresentations of 
those induced from the Borel subgroup ${\rm B}_p$.

\begin{figure}
\centering
\includegraphics[width=0.6\columnwidth]{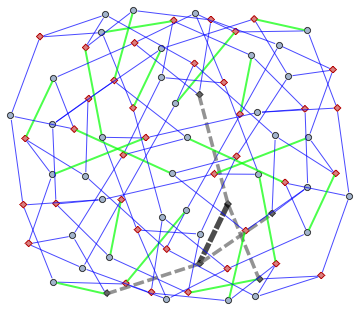}
\caption{The graph of ${\rm SL}(2,\bF_5)$. 
The 120 edges are the group elements, the identity being shown in black.
Its four neighbors (grey) are $-R^{-1}L$, $-L^{-1}R$ and their inverses.
Green: Borel subgroup.
The graph is bipartite, with $V^+$ and $V^-$ colored differently. 
It has the (large) girth = 10 and is Ramanujan.}%
\label{fig:borel}%
\end{figure}
%
\begin{remark} [Principal series reps and the Hilbert space $\tilde{h}(p)$]
 \label{rem:principalseries}
\qquad \newline
The appearance of the principal series representation for the 
Hilbert space $\tilde{h}(p) = \ell^2\big(\Lambda(p)\big)$ is explained as follows.
\begin{enumerate}
\item 
First, we consider $\tilde{h}(n)$ for $n\in \bN$ 
as the subspace of $\ell^2({\rm G}_n)$,
invariant under the right action of the unipotent subgroup 
$U_n$. As such, it is invariant under the 
left action $\hat L_g$ of the $g\in {\rm G}_n$ on $\ell^2({\rm G}_n)$: 
If $\hat R_{\bsm 1&b\\0 & 1\esm}f=f$, then 
\[\hat R_{\bsm 1&b\\0 & 1\esm}\hat L_gf(x)
= (\hat L_gf)\big(x\bsm 1&b\\0 & 1\esm\!\big)
= f\big(g^{-1}x\bsm 1&b\\0 & 1\esm\!\big)
= f(g^{-1}x) = \hat L_gf(x).\]
\item 
For a generator  $\alpha$ of the cyclic group $(\bF_p^\times,\cdot)$
the characters on $\bF_p^\times$ are 
\[\psi_j:\bF_p^\times\to S^1\qmbox{,}\psi_j(\alpha^k)=\exp\big(2\pi\imath jk/(p-1)\big)
\qquad \big(j\in\{0,\ldots,p-2\}\big).\]
Considering  $\ell^2(B_n)$ as the subspace of $\ell^2({\rm G}_n)$ 
given by the functions vanishing outside the Borel subgroup, 
for all $j$ and primes $p$ the character 
\[\tilde \psi_j:B_p\to S^1 \qmbox{,}
\tilde \psi_j\big(\!\bsm a&b\\0 & a^{-1}\esm\!\big) = \psi_j(a)\] 
is in $\ell^2(B_p)\cap \tilde{h}(p)$.
\item 
More generally, the subspaces 
\beq
{\rm Ind}_{p,j} :=
\big\{f\in\ell^2({\rm G}_p)\mid
\forall g\in {\rm G}_p,b\in B_p: f(gb)=\tilde \psi_j(b)f(g) \big\} 
\Leq{indu:rep}
are contained in $\tilde{h}(p)$, and give rise to
the induced representations 
\footnote{Observe the non-standard use of left and right actions.} 
\[\rho_{p,j}:{\rm G}_p\to {\rm GL}({\rm Ind}_{p,j}) \qmbox{,} 
\big(\rho_{p,j}(g)f\big)(x)=f(g^{-1}x)\]
of  the Borel group $B_p$.
\item  
As the $(p+1)$--dimensional subspaces ${\rm Ind}_{p,j}$, $j\in\{0,\ldots,p-2\}$
of $\tilde{h}(p)$ are mutually orthogonal, and  
$\dim\big(\tilde{h}(p)\big)=p^2-1$,
we obtain the orthogonal decomposition
\[\tilde{h}(p) = {\textstyle \bigoplus_{j=0}^{p-2}}\ {\rm Ind}_{p,j}.\hfill \tag*{$\Box$}\]
\end{enumerate}
\end{remark}
The above family of induced representations contains every irreducible 
principal series representation derived from the character $\psi_j$
with multiplicity two if $\psi_j^2\neq 1$ (since $\psi_{p-1-j}=\psi_j^{-1}$), 
and with multiplicity one otherwise.
\subsubsection{Representations of 
${\rm SL}(2,\bZ/n\bZ)$ for Prime Powers {\large $n=p^k$}}
\label{Prime:Powers}
In Nobs and Wolfart \cite{NWI,NWII} all irreducible representations 
for $n=p^k$, $p\in \bP$, are determined.
The method (going back to Weil,  to Kloosterman and to Tanaka) 
to find the representations is to consider 
transformation properties of theta functions associated to binary quadratic forms.
The case $p=2$ needs special treatment. This is an issue to be clarified, 
since \cite{BV}, see below, use in their Section 2 \ Lemma 7.1 of \cite{BG2}, 
which only refers to odd prime powers. 

Inspection of the tables in Section 9 of \cite{NWII}, however, shows that
the dimensions of all such {\em new}  representations 
(the old ones being those already arising for $p^\ell$, $0\le\ell<k$) 
are bounded below by $3n/16$, {\em including} the case $n=2^k$.  
Although \cite{BV} use a lower bound $n/3$ (valid for the odd prime powers), 
this does not change the argument,
since linear growth in $n$ implies
the existence of an $n$--independent spectral gap.
\subsubsection{Representations of 
${\rm SL}(2,\bZ/n\bZ)$ for Arbitrary Integers {\large $n$}}
\label{arbitrary:n}
Using (\ref{congo}), ${\rm G}_n$ is isomorphic to a direct product of the 
${\rm G}_q$, with
the prime power $q$ appearing in the factorization of $n$.
So the irreducible (unitary) representations of ${\rm G}_n$
arise as tensor products of irreducible representations. 

As a consequence of (\ref{dim:SL}) and Sect.\ \ref{Prime:Powers}, the dimensions of
faithful representations are bounded below by $c_\alpha n^\alpha$, for any 
$\alpha\in(0,1)$. This follows from the estimate $\omega(n)=o(\log(n))$ 
for the number of distinct prime factors of $n$, and can
be used in going from the case of prime powers $n$ to
general $n\in\bN$, see the proof of Theorem 1 in Bourgain and Varj\'u \cite{BV}.

A spectral estimate going back to Sarnak and Xue and used in Bourgain and
Gamburd \cite{BG} uses that 
each irreducible representation appears in the regular representation
with multiplicity its dimension.
\cite{BV} then show, given a 
finite symmetric generating set $S\subseteq {\rm SL}(d,\bZ)$ for its image 
$S_n\subseteq {\rm SL}(d,\bZ/n\bZ)$
the expansion property, uniform in $n\in \bN$.
For our context, this implies that the vertex Laplacians 
$\Delta_{V(n)}$ have a spectral gap uniform in $n\in\bN$.
\begin{prop} \label{prop:finite}
There exists an $\vep>0$, such that for 
${\rm G}_n={\rm SL}(2,\bZ/n\bZ)$
\[{\rm spec}\big(\hT^+_{{\rm G}_n}\big)\subseteq C\cup I_\vep\cup \{\pm\eh, 1\}
\mbox{ and }\
{\rm spec}\big(\tilde {t}^+(n)\big)\subseteq C\cup I_\vep\cup \{\pm\eh, 1\}
\quad(n\in\bN), \]
with the circle $C=S^1/\sqrt{2}$ and 
$I_\vep:=\big\{x\in\bR\mid |x|\in [\eh+\vep,1-\vep]\setminus \{1/\sqrt{2}\}\big\}$.
\end{prop}
{\bf Proof.} 
$\bullet$
The first inclusion follows from the result on $\Delta_{V(n)}$ via the relation 
$\Delta_{E(n)}=3(\idty_{\ell^2({\rm G}_n)^+}-\hB^+_{{\rm G}_n})$
and (\ref{rep:T}).\\
$\bullet$
By considering the Hilbert space $\tilde {h}^+(n)$ as a subspace of 
$\ell^2({\rm G}_n)$ as in Remark~\ref{rem:principalseries}, 
\[\tilde {t}^+(n)= \hT^+_{{\rm G}_n}\rstr_{\tilde {h}^+(n)}\qquad(n\in \bN).\]
So the statement about the spectra of the operators 
$\tilde {t}^+(n)$ follows from the one for $\hT^+_{{\rm G}_n}$.
\hfill$\Box$
%
%
\subsection{Finite Adeles}\label{subs:fa}
%
We now consider the set $S= \bP$ of places, 
{\em i.e.}  the ring $\widehat \bZ = \bZ_\bP$ of finite 
integral adeles.
The inverse limit 
\[\bZ_\bP=\varprojlim_{n\in\bN} \bZ/n\bZ\] 
gives rise to natural homomorphisms 
\[\pi_n:\bZ_\bP\ar \bZ/n\bZ \qmbox{and} \pi_n^2:\bZ_\bP^2\ar (\bZ/n\bZ)^2
\qquad (n\in\bN).\]
Functions $f\in h(n)=\ell^2\big((\bZ/n\bZ)^2\big)$ pull back to 
$f\circ\pi_n^2\in \cH_\bP$, the Hilbert space for the finite 
integral adeles.
These are locally constant (or {\em Schwartz--Bruhat}) functions, and with 
our normalization convention (Haar measure having total mass one)
the embeddings $h(n)\hookrightarrow \cH_\bP$ 
are isometric. 

If we thus consider $h(n)$ as a subspace of $\cH_\bP$, 
the latter Hilbert space is a (non--direct) sum of
the $h(n)\ (n\in\bN)$, defined in (\ref{tn:tildetn}). This is implied by the facts that
\begin{enumerate}[$\bullet$]
\item 
every Schwartz--Bruhat function on $\widehat \bZ^2$
is a finite linear combination of characteristic functions 
$\idty_{\tiny\bsm a\\b\esm+n\widehat \bZ^2}$ with 
$n\in\bN$ and $a,b\in\widehat \bZ$ (following {\em e.g.} from 
Lemma 5.4.7 in Deitmar \cite{De}),
\item 
the space of Schwartz--Bruhat functions on $\widehat \bZ^2$ is dense in 
$\cH_\bP=L^2\big(\widehat \bZ^2,m_\bP^2\big)$.
\end{enumerate}
To obtain a direct sum decomposition, we first note that the natural homomorphisms
$\bZ/n\bZ\ar \bZ/m\bZ$ for $m|n$ give rise 
to isometric embeddings $h(m)\hookrightarrow h(n)$.

Then the orthogonal decompositions 
\[h(n)=\bigoplus_{m\in\bN:\, m|n} \tilde h(m)\]
(Lemma 3 of \cite{Kn4}) lead to the orthogonal direct sum
\beq
\cH_\bP=\bigoplus_{m\in\bN}\tilde h(m).
\Leq{nods}
Eigenvectors $f \in h(n)$ of $\htt(n)$ pull back to
eigenvectors $f\circ\pi_n^2\in \cH_\bP$ of $\hT_\bP$, with the same eigenvalue.

We split the Hilbert spaces and operators 
orthogonally into 
\[h(n)= h^+(n)\oplus h^-(n) \qmbox{respectively} 
\htt(n)=\htt^+(n)\oplus \htt^-(n).\]
$h^\pm(n)$ can be considered as finite dimensional subspaces of 
$\cH_\bP^\pm$, and
$\htt^\pm(n)$ the restriction of $\hT_\bP^\pm$.

In analogy with the {\em oldform / newform} calculus of the theory of 
modular forms we need only determine the action of  $\hT_\bP^\pm$ on 
the subspaces $\tilde h^\pm(n):= \tilde h(n)\cap h^\pm(n)$ 
of $h^\pm(n)$.

Then the restricted operators
$\tilde t^\pm(n):=\htt^\pm(n)\rstr_{\tilde h^\pm(n)}$ have the property that
\beq
\bigoplus_{k\in\bN: \; k|n} \tilde h^\pm(k) \cong h^\pm(n)\qmbox{and}
\bigoplus_{k\in\bN: \; k|n} \tilde t^\pm(k) \cong \htt^\pm(n).
\Leq{restricted:op}

So the operators 
$\mathfrak t^\pm(n):= \bigoplus_{k\in\bN: \; k\le n} \tilde t^\pm(k)\ (n\in\bN)$ 
converge strongly to $\hT_\bP^\pm$, and
\beq
{\rm spec}\big(\mathfrak  t^\pm(n)\big) 
\subseteq {\rm spec}\big(\mathfrak  t^\pm(n+1)\big) 
\subseteq {\rm spec}\big(\hT_\bP^\pm\big).
\Leq{spec:inclusion}
Similar to the previous sections, we consider the operators
$\hB_\bP^{+,\pm}$ and $\mathfrak  b^{+,\pm}(n)$ 
$(n\in\bN)$, defined by (\ref{bppm}) and related to $\hT_\bP^+$
respectively $\mathfrak  t^+(n)$ via (\ref{TJB:formula}):
\beq 
\hB_\bP^{+,\pm} = \mp \ed(2\hT_\bP^{+,\pm}+(\hT_\bP^{+,\pm})^{-1})
\mbox{ and } 
\mathfrak  b^{+,\pm}(n) =
\mp \ed \big(2\mathfrak  t^{+,\pm}(n) + (\mathfrak  t^{+,\pm}(n))^{-1}\big).
\Leq{prime:B}
\begin{remark}[Spectra for direct sums] 
The following lemma is not {\em a priori} obvious, since the operator 
$\hT_\bP$ is not normal. {\em E.g.} there is the example of the direct sum 
$\oplus_{k=2}^\infty N_k$ of nilpotent Jordan matrices $N_k$ of size $k$ which
has the closed unit disk as its spectrum (Problem 98 in Halmos, \cite{Ha1}).
\hfill $\Diamond$
\end{remark}
\begin{lemma}
$\overline {\bigcup_{n\in\bN}{\rm spec}\big(\mathfrak  t^+(n)\big)} = 
{\rm spec}\big(\hT_\bP^+\big)$. 
\end{lemma}
{\bf Proof.} 
$\bullet$
Closedness of spectra and (\ref{spec:inclusion}) imply 
$\overline {\bigcup_{n\in\bN}{\rm spec}\big(\mathfrak  t^\pm(n)\big)} \subseteq
{\rm spec}\big(\hT_\bP^\pm\big)$. \\
$\bullet$
By (\ref{TJB:formula}) the operators  
$\mathfrak  b^{+,\pm}(n)$ and 
$\hB_\bP^{+,\pm}$ 
are self-adjoint and bounded in norm by 1.
Holomorphic functional calculus shows that their spectra
have inclusion properties analogous to (\ref{spec:inclusion}), and 
$s-\lim_{n\ar\infty}\mathfrak  b^{+,\pm}(n)= \hB_\bP^{+,\pm}$. So 
by using the resolvent estimate for the normal operators $\mathfrak  b^{+,\pm}(n)$ 
with $\lambda\in\bC$ in the resolvent set
\[\big\|(\mathfrak b^{+,\pm}(n)-\lambda\idty)^{-1}\big\| =
1 / {\rm dist}\big(\lambda,{\rm spec}(\mathfrak b^{+,\pm}(n))\big),\] 
we obtain $\overline {\bigcup_{n\in\bN}{\rm spec}\big(\mathfrak  b^{+,\pm}(n)\big)} = 
{\rm spec}\big(\hB_\bP^{+,\pm}\big)$. (\ref{rep:T}) implies the lemma.
\hfill$\Box$\\[2mm]
%
\begin{figure}[h]%
\centering
\includegraphics[width=0.49\columnwidth]{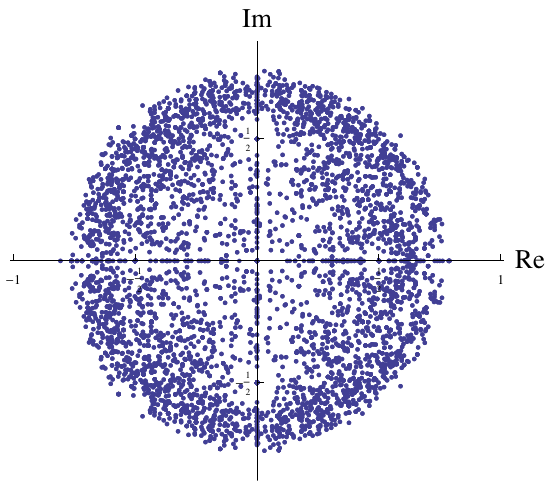} 
\hfill
\includegraphics[width=0.49\columnwidth]{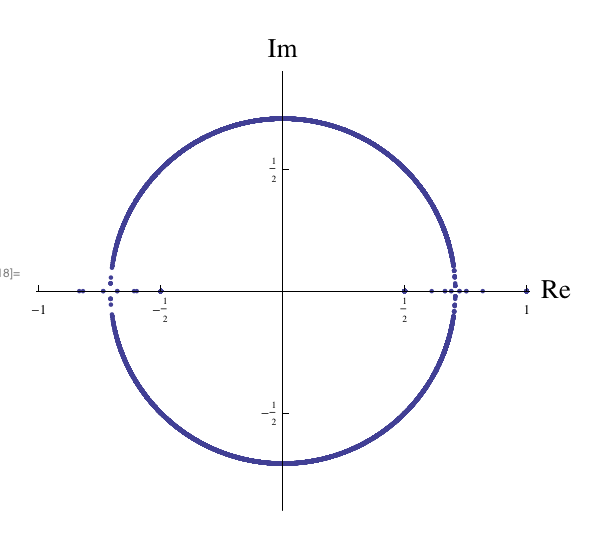} 
\caption{Spectra of the operators $\mathfrak  t^-(50)$ (left) and 
$\mathfrak  t^+(50)$ (right), being subsets of ${\rm spec}\big(\hT_\bP^\pm\big)$ }%
\label{fig:spec}%
\end{figure}
%
Such subsets of the ${\rm spec}\big(\hT_\bP^\pm\big)$ are shown in 
Figure \ref{fig:spec}. They imply for $ \hT_\bP^-$ a spectral radius larger than
$4/5$ (but smaller than $\sqrt{3}/2$ by Proposition \ref{prop:algebraic}), and 
for $ \hT_\bP^+$ (which has an eigenvalue one) the same lower bound for 
the {\em nontrivial} spectral radius, see below.

\begin{prop} \label{nonreal:ok} %
The numbers $\pm\eh$ and $1$ are eigenvalues of infinite multiplicity 
of the operator $\hT_\bP^+$ on $\cH_\bP^+$, and there is a $\vep>0$ with 
\[{\rm spec}\big(\hT_\bP^+\big)\ \subseteq\ 
\{-\eh,\eh,1\}\cup C\cup I_\vep \qmbox{and}
{\rm spec}\big(\hT_\bP^+\big)\cap I_\vep \neq \es\]
with $C=S^1/\sqrt{2}$ and 
$I_\vep=\big\{x\in\bR\mid |x|\in [\eh+\vep,1-\vep]\setminus \{1/\sqrt{2}\}\big\}$.
\end{prop}
{\bf Proof.} 
$\bullet$
Since for all integers $n\in\bN$, $\Lambda(n)\subseteq (\bZ/n\bZ)^2$ 
is a single ${\rm SL} (2,\bZ)$--orbit, and ${\rm SL}(2,\bZ)$ is generated 
by $L$ and $R$, $1$ is an eigenvalue of multiplicity one for the operator 
$\tilde t^+(n)$ on the subspace $\tilde h^+(n)\subseteq \ell^2(\Lambda(n))$. \\
So 1 is of multiplicity $\sigma_0(n)$ for $\htt^+(n)$ (with the divisor function 
$\sigma_0$),
of multiplicity $n$ for $\mathfrak t^+(n)$ and of infinite multiplicity for $\hT_\bP^+$.\\
$\bullet$
That $\pm \eh$ are eigenvalues of $\hT_\bP^+$ can already be inferred from
the case $n=2$ with Hilbert space $h^+(2)=h(2)$. The matrix 
$\htt(2)=\eh \bsm
 2 & 0 &\, 0 & 0 \\
 0 & 1 &\, 0 & 1 \\[1mm]
 0 & 0 &\, 1 & 1 \\
 0 & 1 &\, 1 & 0
\esm$
(w.r.t.\ the basis given by the lexical listing $(00,01,10,11)$ of $(\bZ/2\bZ)^2$)
has the eigenvalues 1, 
$-\eh$ with eigenvector $\bsm0 & -1 &\, -1 & 2\esm^\top$ and 
$\eh$ with eigenvector $\bsm 0 & -1 &\, 1 & 0\esm^\top$. 

In Proposition 11 of \cite{Kn4} the multiplicities of the eigenvalues 
$\pm \eh$ were calculated for
all $n\in\bP$ in 'projective' (that is, dilation invariant) subspaces of $\tilde h(n)$, 
using quadratic reciprocity. In particular these multiplicities are positive for both
signs and all $n\in\bP\setminus\{3,7\}$. \\
This shows that the eigenvalues $\pm \eh$ of
$\hT_\bP^+$ are of infinite multiplicity.\\
$\bullet$
The first {\em nontrivial} (that is, on $I_\vep$) real eigenvalues of 
$h^+(n)$  arise for $n=34$.
They are the real roots of the polynomial 
$64 x^{12}-64 x^{11}+64 x^{10}-64 x^9+36 x^8-26 x^7+18 x^6-13 x^5
+9 x^4-8 x^3+4 x^2-2 x+1$ and equal to 
$0.819427\ldots$ respectively $0.610182\ldots=1/(2\times 0.819427\ldots)$, 
symmetric w.r.t.\ the circle $C$.\\
$\bullet$
That the part  
$\big({\rm spec}\big(\hT_\bP^+\big)\setminus \{\pm\eh, 1\}\big)\cap\bR$ 
of the spectrum is included in $I_\vep$ for some $\vep>0$, is the generalization of 
Proposition 15 of \cite{Kn4} from the projective subspace to the Hilbert space $\tilde h^+(n)$, and from primes to general integers $n$.\\
The argument in \cite{Kn4} used expander estimates for the 
Laplacians $3(\idty-\mathfrak  b^{+,\pm}(n))$, $n\in\bP$, based on
Selberg's Theorem for the congruence subgroups of ${\rm SL}(2,\bZ)$
(see Lubotzky \cite{Lu}, Section 4.3 and 4.4). 
Here we use for $n\in\bN$ Proposition \ref{prop:finite},
based on the estimates by Bourgain, Gamburd and Varj\'u, as explained in
Section \ref{arbitrary:n} above.
\hfill $\Box$
\begin{remark}[Spectral gap] 
Proposition \ref{nonreal:ok} says that $\hT_\bP^+$ has a spectral gap, but
the gap is not of maximal possible size, since 
${\rm spec}\big(\hT_\bP^+\big)\cap I_\vep \neq \es$, 
that is, eigenvalues $\lambda\neq1$ of modulus $|\lambda|>1/\sqrt{2}$
occur. 

One mechanism to reconcile this with RH could be that eigenvalues
of $\htt^+(n)$, which are not equal $\pm\eh$ and do
not already occur for $\htt(m)$ with $m|n$, have modulus
going to $1/\sqrt{2}$ as $n\to\infty$. \hfill $\Diamond$
\end{remark}
By Proposition \ref{nonreal:ok} the operator $\hT_\bP^+$ has a highly 
degenerate spectrum. We now partly lift that degeneracy without changing the
spectrum, by restricting it to the ${\rm SL}(2,\bZ)$--invariant subspace
\[\cH_{\Lambda_\bP} := L^2\big(\Lambda_\bP, m_\bP^2\rstr_{\Lambda_\bP}\big)
\qmbox{of} \cH_\bP = L^2\big(\widehat \bZ^2,m_\bP^2\big).\]
Here 
\beq
\Lambda_\bP:=\prod_{p\in\bP} \Lambda_p\qmbox{with}
\Lambda_p:=\bZ_p^2\setminus p\,\bZ_p^2
\qmbox{equals }\varprojlim_{n\in\bN} \Lambda(n),
\Leq{LambdaP}
the inverse limit being defined w.r.t.\ the homomorphisms 
\[\pi_{n,m}:(\bZ/n\bZ)^2\to (\bZ/m\bZ)^2\qmbox{,}
\bsm\ell+n\bZ\\r+n\bZ\esm\mapsto\bsm\ell+m\bZ\\r+m\bZ\esm 
\qquad(m|n\in\bN), \]
restricted to $ \Lambda(n)$.
By the product formulas \eqref{LambdaP} for $\Lambda_\bP$
and \eqref{def:mS} for $m_\bP$ the restricted Haar
measure $m_\bP^2\rstr_{\Lambda_\bP}$ has total mass $1/\zeta(2)$.
By ${\rm SL}(2,\bZ)$--invariance of $\cH_{\Lambda_\bP}$ the operator
\[\hT_{\Lambda_\bP}:= \hT_\bP\rstr_{\cH_{\Lambda_\bP}}\]
acts on the Hilbert space $\cH_{\Lambda_\bP}$.

It is related to the (non-discrete) Markov chain with state space $\Lambda_\bP$ 
and stochastic kernel
\[\kappa:\Lambda_\bP\times \cB(\Lambda_\bP)\to [0,1]\qmbox{,}
\kappa(x,A)=\eh(\delta_{L(x)}+\delta_{R(x)})(A).\]
The weak Markov property is standard (see, {\em e.g.} Klenke \cite{Kle},
Theorem 17.11). 

The chain is not irreducible in the sense
of Nummelin \cite{Num}.
However, the ${\rm SL}(2,\bZ)$--action on
$(\Lambda_\bP, m_\bP^2\rstr_{\Lambda_\bP})$ is {'irreducible-aperiodic'}
in the weak sense that its projections to $\Lambda(n)$ ($n\in\bN$) are, 
by Lemma 4 of \cite{Kn4}.

Grigorchuk showed in \cite{Gr}, Theorem 1 an individual ergodic theorem for
the Ces\`{a}ro means. In our context it says that for $f\in L^p(\Lambda_\bP)$
with $p\in[1,\infty)$
\[\bar{f}:=\lim_{n\to\infty}\frac{1}{n}\sum_{i=0}^{n-1} 
\hT_{\Lambda_\bP}^i f\ \in \ L^p(\Lambda_\bP),\]
with the same  $L^1(\Lambda_\bP)$ expectation as $f$, and 
$\bar{f}$ is invariant under the action of the
semigroup generated by $\hL$ and $\hR$ on  
$L^p(\Lambda_\bP)$. 
Thus $\bar{f}$ is constant $m_\bP^2\rstr_{\Lambda_\bP}$--a.e.
\begin{lemma}  
${\rm spec}\big(\hT_{\Lambda_\bP}^\pm \big) 
= {\rm spec}\big(\hT_{\bP}^\pm \big)$.\\
The multiplicity of the Perron-Frobenius eigenvalue 1 of $\hT_{\Lambda_\bP}^+$
is one.
\end{lemma}
\textbf{Proof.}
$\bullet$
We only need to show that 
${\rm spec}\big(\hT_{\Lambda_\bP}^\pm \big) \supseteq 
{\rm spec}\big(\hT_{\bP}^\pm\big)$.
For that we lift the eigenfunctions of the $\tilde{t}^\pm(m)$ in
the decomposition (\ref{restricted:op}) to $\tilde{t}^\pm(n)$, 
using~$\pi_{n,m}$.\\
$\bullet$
Clearly the function ${\Lambda_\bP}\to\{1\}$ is eigenfunction of 
$\hT_{\Lambda_\bP}^+$ with eigenvalue 1. 
If, on the other hand,  $f:{\Lambda_\bP}\to\bC$ is an eigenfunction of
$\hT_{\Lambda_\bP}^+$ with eigenvalue 1, then it is left-and right-invariant:
$\hL_{\Lambda_\bP}^+f=f$ and $\hR_{\Lambda_\bP}^+f=f$. So it is
${\rm SL}(2,\bZ)$--invariant, and by the above constant a.e.
\hfill$\Box$\\[2mm]
This restricted operator is important, since $\Lambda_\bP$ is the closure
of the ${\rm SL}(2,\bZ)$--orbit of our initial point
$\bsm1\\1\esm\in \widehat \bZ^2$.
%
\subsection{The Adelic Case}\label{subs:ac}
%
We now consider for $\bP_\infty=\{\infty\}\cup\bP$  
the adelic  Markov operator $\hT_{\bP_\infty}$ on
the Hilbert space
\[\cH_{\bP_\infty}=L^2\big(\bZ^2\times \widehat\bZ^2,m_{\bP_\infty}^2\big).\]
The following statement is the main result of this article. It characterizes
the spectrum of $ \hT_{\bP_\infty}$ as the union 
of spectra of the operator $\hT_{\bP}$
analyzed in Proposition \ref{nonreal:ok} and of an operator 
$\hT_{\Lambda\times \widehat \bZ^2}$, which in turn has the same spectrum as
$\hT_{\Lambda}$ from  Proposition  \ref{prop:real}.
\begin{theorem} \label{theorem:adelic} %
${\rm spec}\big( \hT_{\bP_\infty}^\pm \big)
= {\rm spec}\big( \hT_{\bP}^\pm \big) \cup
{\rm spec}\big(\hT_{\Lambda\times \widehat \bZ^2}^\pm\big)$ with
\[ {\rm spec}\big(\hT_{\Lambda\times \widehat \bZ^2}^+\big) = 
\bigcup_{n\in\bN}{\rm spec}\big(\hT_{\Lambda\times \Lambda(n)}^+\big)
= \{-\eh,\eh\}\cup C.\]
In particular the adelic Markov operator $\hT_{\bP_\infty}^+$ has a spectral gap.
\end{theorem}
{\bf Proof.} 
$\bullet$
Because of (\ref{split:real}) and  (\ref{nods}) the Hilbert space
\[\cH_{\bP_\infty}=L^2\big(\bZ^2\times \widehat\bZ^2,m_{\bP_\infty}^2\big)
\cong \ell^2(\bZ^2)\otimes  L^2\big(\widehat\bZ^2,m_\bP^2\big)\]
splits into a non-orthogonal \marginpar{} direct sum 
\[\cH_{\bP_\infty}\cong\bigoplus_{m\in\bN_0,\, n\in\bN}\ell^2(m\,\Lambda)\;\otimes\; 
\ell^2\big(\Lambda(n)\big)
\cong\bigoplus_{m\in\bN_0,\, n\in\bN}\ell^2\big(m\Lambda\times
\Lambda(n)\big).\]
The operator $\hT_{\bP_\infty}$ splits accordingly into a direct sum of 
$\hT_{m\Lambda\times\Lambda(n)}$.

Since $\ell^2(0\,\Lambda)\cong\bC$ and $\ell^2(m\,\Lambda)\cong\ell^2(\Lambda)$
(see the proof of Proposition \ref{prop:real})
\beq
\hT_{0\Lambda\times\Lambda(n)}\cong\htt(n) \qmbox{and}
\hT_{m\Lambda\times\Lambda(n)}\cong\hT_{\Lambda\times\Lambda(n)}
\quad(m,n\in\bN).
\Leq{decompose}
$\bullet$
The first identity follows from the direct sum decomposition of 
$\hT_{\bP_\infty}^\pm $ and (\ref{decompose}).\\
$\bullet$
$\hT_{\Lambda\times \Lambda(n)}^+
=\eh(\hL_{\Lambda\times \Lambda(n)}^+ +\hR_{\Lambda\times \Lambda(n)}^+)
\cong\eh(\hL_{\Lambda}^+ \otimes \hL_{\Lambda(n)}^+ 
+\hR_{\Lambda}^+ \otimes \hR_{\Lambda(n)}^+)$ 
has for all $n\in\bN$ the eigenvalues $\pm\eh$, since according to 
Prop.\ \ref{prop:real}
$\hT_{\Lambda}^+=\eh(\hL_{\Lambda}^+ + \hR_{\Lambda}^+)$ 
has such eigenfunctions $\phi_\pm$, and for the constant eigenfunction
$\idty_{\Lambda(n)}\in h(n)$ of $\hL_{\Lambda(n)}^+$ and
$\hR_{\Lambda(n)}^+$ with eigenvalue one $\phi_\pm\otimes\idty_{\Lambda(n)}$
are eigenfunctions of $\hT_{\Lambda\times \Lambda(n)}^+$ with eigenvalues 
$\pm\eh$.\\
$\bullet$
By a converse argument $\hT_{\Lambda\times \Lambda(n)}^+$ does not have the
eigenvalues $\pm1$, since these would imply that $\hT_{\Lambda}^+$ had these
eigenvalues, contradicting Proposition \ref{prop:real}.\\
$\bullet$
To show that ${\rm spec}\big(\hT_{\Lambda\times \Lambda(n)}^+\big)
\subseteq\{-\eh,\eh\}\cup C$, we first consider the resolvent of 
$\hB_{{\rm SL}\times \Lambda(n)}^+$. Arguing in a way analogous to
the one in the proof of Prop.\ \ref{prop:SL}, we obtain the estimate
\beq
\|D_{{\rm SL}\times \Lambda(n)}(k)\| \le 3\,n\,k\, 2^{k/2}\qquad(k\in\bN)
\Leq{Dk:est:adelic}
for the operators (indexed by $k\in\bN_0$)
\[D(k)\equiv D_{{\rm SL}\times \Lambda(n)}(k)\in 
\hB(\cH_{V,{\rm SL}\times \Lambda(n)})
\qmbox{,}
\big(D(k)f\big)(v):= \sum_{w:\;{\rm dist}(v,w)=k} f(w)\]
on the vertex Hilbert space $\cH_{V,{\rm SL}\times \Lambda(n)}$ of
the $\hY_\pm$ action on ${\rm SL}(2,\bZ)\times \Lambda(n)$.
The constant factor $n$ in the estimate (\ref{Dk:est:adelic}) of the 
$D_{{\rm SL}\times \Lambda(n)}(k)$ does not change the
convergence properties of the resolvent of 
$\hB_{{\rm SL}\times \Lambda(n)}^+$ on $\cN^+$, compared to the one for 
$\hB_{{\rm SL}}^+$. So ${\rm spec}\big( \hB_{{\rm SL}\times \Lambda(n)}^+\big)
\subseteq [-\sqrt{8}/3,\sqrt{8}/3]$. 
We employ the commutative diagram
\[\begin{CD}
U @>\hB_{{\rm SL}\times \Lambda(n)}>>U\\
@A\widehat E AA@VV\widehat  \Pi V\\
\ell^2\big(\Lambda\times \Lambda(n)\big) @>\hB_{\Lambda\times \Lambda(n)}>> 
\ell^2\big(\Lambda\times \Lambda(n)\big)
\end{CD}\]
defined like (\ref{CD1}) to conclude that 
${\rm spec}\big( \hB_{\Lambda\times \Lambda(n)}^+\big)\subseteq
[-\sqrt{8}/3,\sqrt{8}/3]$, too.
\hfill $\Box$\\[2mm]

\noindent
{\bf Acknowledgement:} I thank Johannes Singer (Erlangen) for his comments
and Gunther Cornelissen (Utrecht) for showing me the references \cite{NWI,NWII}.

%
\addcontentsline{toc}{section}{Literature}
%

\end{document}